\documentclass{article}

\title{Rooted tree graphs and the Butcher group:\\
Combinatorics of elementary perturbation theory}

\author{William G. Faris\\
NYU Shanghai and University of Arizona}

\newtheorem{theorem}{Theorem}
\newtheorem{corollary}{Corollary}
\newtheorem{proposition}{Proposition}
\newenvironment{proof}{\textbf{Proof:}\ }{\ $||$}
\newenvironment{example}{\textbf{Example:}\ }{\ $||$}
\newenvironment{remark}{\textbf{Remark:}\ }{\ $||$}

\newcommand{\residue}{\mathrm{res}}
\newcommand{\leqtree}{\leq_{\mathrm{tree}}}
\newcommand{\leqlin}{\leq_{\mathrm{lin}}}

\newcommand{\negthreegl}{-120}
\newcommand{\negtwogl}{-80}
\newcommand{\neggl}{-40}
\newcommand{\neghalfgl}{-20}
\newcommand{\halfgl}{20}
\newcommand{\gl}{40}
\newcommand{\twogl}{80}

\newcommand{\fourgl}{160}

\newcommand{\sixgl}{240}

\newcommand{\bggl}{280}

\newcommand{\negfoursgl}{-80}

\newcommand{\negtwosgl}{-40}
\newcommand{\negsgl}{-20}
\newcommand{\neghalfsgl}{-10}

\newcommand{\sgl}{20}
\newcommand{\twosgl}{40}
\newcommand{\threesgl}{60}
\newcommand{\foursgl}{80}

\newcommand{\sixsgl}{120}
\newcommand{\sevensgl}{140}
\newcommand{\eightsgl}{160}
\newcommand{\ninesgl}{180}
\newcommand{\tensgl}{200}
\newcommand{\twelvesgl}{240}
\newcommand{\bgsgl}{140}

\newcommand{\setpartition}{\mathrm{Part}}

\newcommand{\Rb}{\mathbf{R}}

\newcommand{\Ac}{\mathcal{A}}
\newcommand{\Bc}{\mathcal{B}}
\newcommand{\Gc}{\mathcal{G}}
\newcommand{\Hc}{\mathcal{H}}

\begin{document}

\maketitle

\begin{abstract}
The perturbation expansion of the solution of a fixed point equation or of an ordinary differential equation may
be expressed as a power series in the perturbation parameter. The terms in this series are indexed by rooted trees and
depend on a parameter in the equation in a way determined by the structure of the tree. Power series of this form may be considered more
generally; there are two interesting and useful group structures on these series, corresponding to operations of composition
and substitution. The composition operation defines the Butcher group, an infinite dimensional group that was first introduced in the context of numerical analysis. This survey discusses
various ways of realizing these rooted trees: as labeled rooted trees, or increasing labeled rooted trees, or
unlabeled rooted trees. It is argued that the simplest framework is to use labeled rooted trees.
\end{abstract}

\section{Introduction}

This year we celebrate the scientific contributions of Charles Newman. Chuck has worked in
almost every aspect of mathematical physics. He can identify a significant problem, locate the appropriate framework,
and find an unexpected path to a comprehensible solution.  His insights
and his generosity in sharing them are extraordinarily valuable to the community.
Chuck has been my friend and colleague from years together  at the University of Arizona and now at NYU Shanghai.
It is a pleasure to dedicate this paper to such a distinguished scientist.

The paper is a largely expository survey of rooted trees and the Butcher group. The Butcher group is an
infinite dimensional group associated with rooted trees.  See \cite{ChartierHairerVilmart}
for the current status of this subject. Most expositions are in the framework of unlabeled rooted trees.
The main message of the present work is that the combinatorics is simpler when
formulated in terms of labeled rooted trees. In particular, group structures associated with rooted trees
arise naturally from calculations using elementary calculus formulas. This is in the spirit of the theory
of combinatorial species \cite{BergeronLabelleLeroux}. This same point of view is fruitful in the study
of graph expansions in statistical mechanics \cite{Faris1} and of diagram expansions in quantum field theory \cite{Faris2}.

The subject matter of the present work has two independent origins. One begins in 1963 with work of Butcher
in numerical analysis. He discovered that a large class of numerical methods for ordinary differential equations
may be expressed as sums indexed by rooted trees. Furthermore, these sums may be combined in a way that
defines a group structure. This subject has become part of the lore of numerical analysis \cite{Butcher,HairerWanner1,HairerWanner2}.
It remains active; for instance see \cite{LundervoldMuntheKaas,BogfjellmoSchmeding}    and works cited therein.

The other origin is research on renormalization in quantum field theory, starting with the 1998
contribution of  Connes and Kreimer \cite{ConnesKreimer}. This work is usually presented in the language
of combinatorial Hopf algebras. Many authors, for instance  \cite{Brouder,GirelliKrajewskiMartinetti},
have investigated the relation between the Butcher group and problems in quantum field theory. Recently
there has been an explosion of mathematics papers  treating  Hopf algebras associated with rooted trees.
See for instance \cite{CalaqueEbrahimiFardManchon} and papers cited there.

The approach here  begins by distinguishing basic rooted tree constructions. The
starting point is a finite set $U$ known as the \emph{label set} or \emph{vertex set}.
A labeled rooted tree on $U$ is a tree graph with vertex set $U$ together with a distinguished
point $r$ in $U$. An increasing rooted tree is one for which the label set $U$ is  linearly ordered and the labels increase with distance from
the root $r$. An unlabeled rooted tree  is an isomorphism class of labeled rooted trees.
 \begin{itemize}
\item $\Ac[U]$ is all  rooted trees on  label set $U$.
\item $\Ac^\uparrow[U]$ is all increasing rooted trees on linearly ordered  label set $U$.
\item $\tilde \Ac[n]$ is  all unlabeled rooted trees with $n$ vertices.
\end{itemize}
The $\Ac$ notation is from the French ``arbre''.
The relation between these constructions is that if $U$ has $n$ elements, then
\begin{equation}
\Ac^\uparrow[U] \to \Ac[U] \to \tilde \Ac[n].
\end{equation}
The first map is an injection, and the second map is a surjection.

The topics discussed include:
\begin{itemize}
\item Unlabeled, labeled, and increasing rooted trees
\item Fixed point equations and ordinary differential equations
\item The composition operation (Butcher group)
\item The substitution operation
\end{itemize}
An appendix reviews  calculus formulas as used in combinatorics.

\begin{example} Figure~1  illustrates the distinction between  labeled rooted trees and unlabeled rooted trees in the case of a vertex set with
three elements. There are 9 labeled rooted trees.   However there are only 2 unlabeled  rooted trees. One of these has a symmetry that exchanges the two
non-root vertices.
\end{example}

% Figure 1 Labeled rooted trees on three vertices
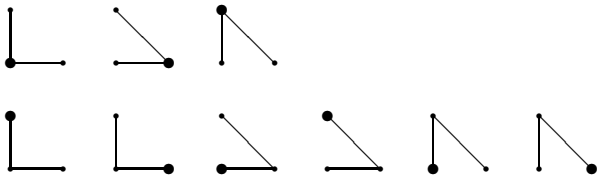
\begin{figure}[tb]
\begin{picture}(\bgsgl,\twosgl)(\negfoursgl,0)
\put(0,\twosgl){
\begin{picture}(\sgl,\sgl)
\put(0,0){\line(1,0){\sgl}}
\put(0,0){\line(0,1){\sgl}}
%\put(0,0){\line(1,1){\sgl}}
%\put(\sgl,\sgl){\line(-1,0){\sgl}}
%\put(\sgl,\sgl){\line(0,-1){\sgl}}
%\put(0,\sgl){\line(1,-1){\sgl}}
\put(0,0){\circle*{4}}
\put(\sgl,0){\circle*{2}}
\put(0,\sgl){\circle*{2}}
%\put(\sgl,\sgl){\circle*{2}}
%\put(\neghalfsgl,\halfsgl){4}
\end{picture}}
\put(\twosgl,\twosgl){
\begin{picture}(0,0)
\put(0,0){\line(1,0){\sgl}}
%\put(0,0){\line(0,1){\sgl}}
%\put(0,0){\line(1,1){\sgl}}
%\put(\sgl,\sgl){\line(-1,0){\sgl}}
%\put(\sgl,\sgl){\line(0,-1){\sgl}}
\put(0,\sgl){\line(1,-1){\sgl}}
\put(0,0){\circle*{2}}
\put(\sgl,0){\circle*{4}}
\put(0,\sgl){\circle*{2}}
%\put(\sgl,\sgl){\circle*{2}}
%\put(\neghalfsgl,\halfsgl){24}
\end{picture}}
\put(\foursgl,\twosgl){
\begin{picture}(\sgl,\sgl)
%\put(0,0){\line(1,0){\sgl}}
\put(0,0){\line(0,1){\sgl}}
%\put(0,0){\line(1,1){\sgl}}
%\put(\sgl,\sgl){\line(-1,0){\sgl}}
%\put(\sgl,\sgl){\line(0,-1){\sgl}}
\put(0,\sgl){\line(1,-1){\sgl}}
\put(0,0){\circle*{2}}
\put(\sgl,0){\circle*{2}}
\put(0,\sgl){\circle*{4}}
%\put(\sgl,\sgl){\circle*{2}}
%\put(\neghalfsgl,\halfsgl){12}
\end{picture}}
\put(0,0){
\begin{picture}(\sgl,\sgl)
\put(0,0){\line(1,0){\sgl}}
\put(0,0){\line(0,1){\sgl}}
%\put(0,0){\line(1,1){\sgl}}
%\put(\sgl,\sgl){\line(-1,0){\sgl}}
%\put(\sgl,\sgl){\line(0,-1){\sgl}}
%\put(0,\sgl){\line(1,-1){\sgl}}
\put(0,0){\circle*{2}}
\put(\sgl,0){\circle*{2}}
\put(0,\sgl){\circle*{4}}
%\put(\sgl,\sgl){\circle*{2}}
%\put(\neghalfsgl,\halfsgl){24}
\end{picture}}
\put(\twosgl,0){
\begin{picture}(\sgl,\sgl)
\put(0,0){\line(1,0){\sgl}}
\put(0,0){\line(0,1){\sgl}}
%\put(0,0){\line(1,1){\sgl}}
%\put(\sgl,\sgl){\line(-1,0){\sgl}}
%\put(\sgl,\sgl){\line(0,-1){\sgl}}
%\put(0,\sgl){\line(1,-1){\sgl}}
\put(0,0){\circle*{2}}
\put(\sgl,0){\circle*{4}}
\put(0,\sgl){\circle*{2}}
%\put(\sgl,\sgl){\circle*{2}}
%\put(\neghalfsgl,\halfsgl){4}
\end{picture}}
\put(\foursgl,0){
\begin{picture}(0,0)
\put(0,0){\line(1,0){\sgl}}
%\put(0,0){\line(0,1){\sgl}}
%\put(0,0){\line(1,1){\sgl}}
%\put(\sgl,\sgl){\line(-1,0){\sgl}}
%\put(\sgl,\sgl){\line(0,-1){\sgl}}
\put(0,\sgl){\line(1,-1){\sgl}}
\put(0,0){\circle*{4}}
\put(\sgl,0){\circle*{2}}
\put(0,\sgl){\circle*{2}}
%\put(\sgl,\sgl){\circle*{2}}
%\put(\neghalfsgl,\halfsgl){24}
\end{picture}}
\put(\sixsgl,0){
\begin{picture}(\sgl,\sgl)
\put(0,0){\line(1,0){\sgl}}
%\put(0,0){\line(0,1){\sgl}}
%\put(0,0){\line(1,1){\sgl}}
%\put(\sgl,\sgl){\line(-1,0){\sgl}}
%\put(\sgl,\sgl){\line(0,-1){\sgl}}
\put(0,\sgl){\line(1,-1){\sgl}}
\put(0,0){\circle*{2}}
\put(\sgl,0){\circle*{2}}
\put(0,\sgl){\circle*{4}}
%\put(\sgl,\sgl){\circle*{2}}
%\put(\neghalfsgl,\halfsgl){12}
\end{picture}}
\put(\eightsgl,0){
\begin{picture}(\sgl,\sgl)
%\put(0,0){\line(1,0){\sgl}}
\put(0,0){\line(0,1){\sgl}}
%\put(0,0){\line(1,1){\sgl}}
%\put(\sgl,\sgl){\line(-1,0){\sgl}}
%\put(\sgl,\sgl){\line(0,-1){\sgl}}
\put(0,\sgl){\line(1,-1){\sgl}}
\put(0,0){\circle*{4}}
\put(\sgl,0){\circle*{2}}
\put(0,\sgl){\circle*{2}}
%\put(\sgl,\sgl){\circle*{2}}
%\put(\neghalfsgl,\halfsgl){4}
\end{picture}}
\put(\tensgl,0){
\begin{picture}(0,0)
%\put(0,0){\line(1,0){\sgl}}
\put(0,0){\line(0,1){\sgl}}
%\put(0,0){\line(1,1){\sgl}}
%\put(\sgl,\sgl){\line(-1,0){\sgl}}
%\put(\sgl,\sgl){\line(0,-1){\sgl}}
\put(0,\sgl){\line(1,-1){\sgl}}
\put(0,0){\circle*{2}}
\put(\sgl,0){\circle*{4}}
\put(0,\sgl){\circle*{2}}
%\put(\sgl,\sgl){\circle*{2}}
%\put(\neghalfsgl,\halfsgl){24}
\end{picture}}
\end{picture}
\caption{Labeled rooted trees on three vertices}
\end{figure}

\section{Labeled rooted trees}

\noindent\textbf{Labeled rooted trees as functions}
Let $U$ be a non-empty finite set. Let $f: U \to U$ be a function with a single fixed point $r$ such that $U \setminus \{r\}$ has
no non-empty invariant subset. Then $f$ defines a labeled rooted tree. In the following it will be more convenient to
consider the function  $f$ restricted to $U \setminus \{r\}$. This leads to the official definition of labeled rooted tree  used
here.

Let $U$ be a non-empty finite set. Let $r$ be a point in $U$, and let $T: U \setminus \{r\} \to U$ be a function that has no non-empty invariant subset.
Then $T$ is called a \emph{labeled rooted tree}. The set $U$ is the label set, and the point $r$ is the \emph{root}. A point in $U$ is called a \emph{label}
or a \emph{vertex}. Each ordered pair $(i, T(i))$ with $i \neq r$ is called an \emph{edge}. The point $T(i)$
is called the \emph{predecessor} of $i$.

 The set of labeled rooted trees with label set $U$ is $\Ac[U]$. If $U$ is empty, then there are no labeled rooted trees on $U$, so
 $\Ac[\emptyset] = \emptyset$. If $T$ is in $\Ac[U]$, then $[T] = U$ is the
label set of $T$. The number of points in the label set is $|T| = |U|$.  The set of \emph{immediate successor} points that map to $j$ in $U$ is
$T^{-1}(j)$. The \emph{degree} of $j$ is $|T^{-1}(j)|$, the number of immediate successor points.
(The definition of degree used here is special to rooted trees; it is not the usual definition from graph theory.)
A \emph{leaf} is a vertex with degree zero. For a one-vertex tree $\bullet$ the root is a leaf.

%Consider a labeled rooted tree $T$. For each vertex $j$ there is a rooted tree $T_j$ whose vertex set consists of all vertices that
%are sent to $j$ by some iterate of $T$. The tree $T_j$ is the restriction of $T$ to this set. Thus there are two notions of the size
%of a vertex $j$, the degree $|T^{-1}(j)|$ and the number of points $|T_j|$ in the tree $T_j$ over $j$. Both will turn out to be useful.

If $b: W \to U$ is a bijection, then the map $T \mapsto T \circ b$ maps $\Ac[U]$ to $\Ac[W]$. Such a map is
called a \emph{relabeling}. Most interesting properties of labeled rooted trees are not affected by relabeling.
It might seem reasonable to use a standard label set $U_n$ for each $n$. An obvious candidate is $[1,n]$, that is,
the set $\{1, \ldots, n\}$. On the other hand, it is common to consider a labeled rooted tree on a subset of $U$, and
that means that other labels sets are going to arise naturally.

The set $\Ac$ of all labeled rooted trees may be defined by choosing a label set $U$ with $|U| = n$ for each $n = 1,2,3,\ldots$.
For some purposes it is useful to adjoin an empty set object associated with the empty label set. The set of all labeled
rooted trees with this extra object is written $\Ac_\emptyset$.

\noindent\textbf{Labeled rooted trees as partially ordered sets}
A labeled rooted tree $T$ on a  finite label set $U \neq \emptyset$ may be viewed as a partial order $\leqtree$ on $U$. This is the
unique partial order with the property that $T(j) \leqtree j$ for every vertex $j \neq r$.

For each vertex $i$ there is a rooted tree $T_i$ whose vertex set consists of all vertices that
are sent to $i$ by some iterate of $T$. The tree $T_i$ is the restriction of $T$ to this set.
To say that $ i \leqtree j$ is the same as saying that $j$ is in the vertex set of the tree $T_i$.

The special feature of this partial order is that for each $j$ in $U$ the
set of all $i \leqtree j$ is linearly ordered with respect to the restriction of $\leqtree$ to this set. Furthermore, there is a least
element in $U$, the root $r$. There can be one or more maximal elements; these are the leaves.

\noindent\textbf{Labeled rooted trees as graphs}
A  {labeled tree} $T$ on a non-empty finite set $U$ is a simple graph with $U$ as vertex set that is connected and has no cycles.
For every pair of vertices  $i \neq j$ there is a unique simple path connecting the two points.
A {labeled rooted tree} $T$ on $U$ is equivalent to a labeled tree and a choice of root point in $U$.   For each vertex $i$ other than the root,
there is a unique edge from $i$ in the direction of the root and corresponding vertex $T(i)$. The usual way of picturing a labeled rooted
tree is as a set together with tree graph and distinguished point.

\noindent\textbf{Forests of labeled rooted trees}
Let $V$ be a  finite set. Let $f: V \to V$ be a function with a set of fixed points $R$ such that $U \setminus R$ has
no non-empty invariant subset.
Then consider the function $f$ restricted to $U \setminus R$. This motivates the official definition of forest of labeled rooted trees.

A \emph{forest of labeled rooted trees} on a set $V$ is a subset $R \subseteq V$ and a function $F: V \setminus R \to V$ such that $V \setminus R$
has no non-empty invariant subset. It is possible that $V = \emptyset$, in which case $R = \emptyset$, and $F$ is the empty function.

The most important fact about a forest of labeled rooted trees is that there
is a set partition $\Gamma$ of $V$ with the following property. For each block $B$ of $\Gamma$ there is a unique $r$ in $R$
such that the restriction $F_B$ of $F$ to $B \setminus \{r\}$ is a labeled rooted tree. When $V = \emptyset$ this is the empty
set partition.

\iffalse
Let $R$ be the set of roots of the rooted trees in the forest. A forest on a set $V$ may also be considered as a function defined on
$U \setminus R$. The value of the function on a vertex $i$ in block $B$ is obtained by applying the rooted tree function $F(B)$ to $i$.

Remark: Let $U$ be a finite set. Let  $T: U \to U$ be an arbitrary function from the set to itself. Then there is a set partition $\Gamma$ of $U$
with the following property. For each block $B$ of the set partition there is a subset $C_B \subseteq B$ such that $f$ restricted to $C_B$
is a cycle with some period (possibly  one). The block $B$ itself admits a set partition into blocks $U_j$ for $j$ in $C_B$ so that $T:
U_j \setminus C_B \to U_j$ is a rooted tree function.  In conclusion, an arbitrary function consists of trees feeding into cycles.

Suppose that the function $T:U \to U$ has no non-trivial cycles, only fixed points. Then there is a set partition $\Gamma$ of $U$
with the following property. For each block $B$ of the set partition there is a point $j$ in $ B$ that is a fixed point of $T$.
Furthermore $T: U_j \setminus \{j\} \to U_j$ is a rooted tree function.  Thus an  arbitrary function with no
cycles consists of trees feeding into fixed points. A function with no cycles is essentially the same thing as a forest function.
\fi

% Figure 2 Labeled rooted trees on three vertices: Prufer sequences
\begin{figure}[tb]
\begin{picture}(\bgsgl,\sevensgl)(\negfoursgl,\negsgl)
\put(0,\foursgl){
\begin{picture}(\sgl,\sgl)
\put(0,0){\line(1,0){\sgl}}
\put(0,0){\line(0,1){\sgl}}
%\put(0,0){\line(1,1){\sgl}}
%\put(\sgl,\sgl){\line(-1,0){\sgl}}
%\put(\sgl,\sgl){\line(0,-1){\sgl}}
%\put(0,\sgl){\line(1,-1){\sgl}}
\put(-8,0){1}
\put(23,0){2}
\put(-8,20){3}
%\put(23,\gl){4}
\put(0,0){\circle*{4}}
\put(\sgl,0){\circle*{2}}
\put(0,\sgl){\circle*{2}}
%\put(\sgl,\sgl){\circle*{2}}
%\put(\neghalfsgl,\halfsgl){4}
\put(0,\negsgl){11}
\end{picture}}
\put(\twosgl,\foursgl){
\begin{picture}(0,0)
\put(0,0){\line(1,0){\sgl}}
%\put(0,0){\line(0,1){\sgl}}
%\put(0,0){\line(1,1){\sgl}}
%\put(\sgl,\sgl){\line(-1,0){\sgl}}
%\put(\sgl,\sgl){\line(0,-1){\sgl}}
\put(0,\sgl){\line(1,-1){\sgl}}
\put(-8,0){1}
\put(23,0){2}
\put(-8,20){3}
%\put(23,\gl){4}
\put(0,0){\circle*{2}}
\put(\sgl,0){\circle*{4}}
\put(0,\sgl){\circle*{2}}
%\put(\sgl,\sgl){\circle*{2}}
%\put(\neghalfsgl,\halfsgl){24}
\put(0,\negsgl){22}
\end{picture}}
\put(\foursgl,\foursgl){
\begin{picture}(\sgl,\sgl)
%\put(0,0){\line(1,0){\sgl}}
\put(0,0){\line(0,1){\sgl}}
%\put(0,0){\line(1,1){\sgl}}
%\put(\sgl,\sgl){\line(-1,0){\sgl}}
%\put(\sgl,\sgl){\line(0,-1){\sgl}}
\put(0,\sgl){\line(1,-1){\sgl}}
\put(-8,0){1}
\put(23,0){2}
\put(-8,20){3}
%\put(23,\gl){4}
\put(0,0){\circle*{2}}
\put(\sgl,0){\circle*{2}}
\put(0,\sgl){\circle*{4}}
%\put(\sgl,\sgl){\circle*{2}}
%\put(\neghalfsgl,\halfsgl){12}
\put(0,\negsgl){33}
\end{picture}}
\put(0,0){
\begin{picture}(\sgl,\sgl)
\put(0,0){\line(1,0){\sgl}}
\put(0,0){\line(0,1){\sgl}}
%\put(0,0){\line(1,1){\sgl}}
%\put(\sgl,\sgl){\line(-1,0){\sgl}}
%\put(\sgl,\sgl){\line(0,-1){\sgl}}
%\put(0,\sgl){\line(1,-1){\sgl}}
\put(-8,0){1}
\put(23,0){2}
\put(-8,20){3}
%\put(23,\gl){4}
\put(0,0){\circle*{2}}
\put(\sgl,0){\circle*{2}}
\put(0,\sgl){\circle*{4}}
%\put(\sgl,\sgl){\circle*{2}}
%\put(\neghalfsgl,\halfsgl){24}
\put(0,\negsgl){13}
\end{picture}}
\put(\twosgl,0){
\begin{picture}(\sgl,\sgl)
\put(0,0){\line(1,0){\sgl}}
\put(0,0){\line(0,1){\sgl}}
%\put(0,0){\line(1,1){\sgl}}
%\put(\sgl,\sgl){\line(-1,0){\sgl}}
%\put(\sgl,\sgl){\line(0,-1){\sgl}}
%\put(0,\sgl){\line(1,-1){\sgl}}
\put(-8,0){1}
\put(23,0){2}
\put(-8,20){3}
%\put(23,\gl){4}
\put(0,0){\circle*{2}}
\put(\sgl,0){\circle*{4}}
\put(0,\sgl){\circle*{2}}
%\put(\sgl,\sgl){\circle*{2}}
%\put(\neghalfsgl,\halfsgl){4}
\put(0,\negsgl){12}
\end{picture}}
\put(\foursgl,0){
\begin{picture}(0,0)
\put(0,0){\line(1,0){\sgl}}
%\put(0,0){\line(0,1){\sgl}}
%\put(0,0){\line(1,1){\sgl}}
%\put(\sgl,\sgl){\line(-1,0){\sgl}}
%\put(\sgl,\sgl){\line(0,-1){\sgl}}
\put(0,\sgl){\line(1,-1){\sgl}}
\put(-8,0){1}
\put(23,0){2}
\put(-8,20){3}
%\put(23,\gl){4}
\put(0,0){\circle*{4}}
\put(\sgl,0){\circle*{2}}
\put(0,\sgl){\circle*{2}}
%\put(\sgl,\sgl){\circle*{2}}
%\put(\neghalfsgl,\halfsgl){24}
\put(0,\negsgl){21}
\end{picture}}
\put(\sixsgl,0){
\begin{picture}(\sgl,\sgl)
\put(0,0){\line(1,0){\sgl}}
%\put(0,0){\line(0,1){\sgl}}
%\put(0,0){\line(1,1){\sgl}}
%\put(\sgl,\sgl){\line(-1,0){\sgl}}
%\put(\sgl,\sgl){\line(0,-1){\sgl}}
\put(0,\sgl){\line(1,-1){\sgl}}
\put(-8,0){1}
\put(23,0){2}
\put(-8,20){3}
%\put(23,\gl){4}
\put(0,0){\circle*{2}}
\put(\sgl,0){\circle*{2}}
\put(0,\sgl){\circle*{4}}
%\put(\sgl,\sgl){\circle*{2}}
%\put(\neghalfsgl,\halfsgl){12}
\put(0,\negsgl){23}
\end{picture}}
\put(\eightsgl,0){
\begin{picture}(\sgl,\sgl)
%\put(0,0){\line(1,0){\sgl}}
\put(0,0){\line(0,1){\sgl}}
%\put(0,0){\line(1,1){\sgl}}
%\put(\sgl,\sgl){\line(-1,0){\sgl}}
%\put(\sgl,\sgl){\line(0,-1){\sgl}}
\put(0,\sgl){\line(1,-1){\sgl}}
\put(-8,0){1}
\put(23,0){2}
\put(-8,20){3}
%\put(23,\gl){4}
\put(0,0){\circle*{4}}
\put(\sgl,0){\circle*{2}}
\put(0,\sgl){\circle*{2}}
%\put(\sgl,\sgl){\circle*{2}}
%\put(\neghalfsgl,\halfsgl){4}
\put(0,\negsgl){31}
\end{picture}}
\put(\tensgl,0){
\begin{picture}(0,0)
%\put(0,0){\line(1,0){\sgl}}
\put(0,0){\line(0,1){\sgl}}
%\put(0,0){\line(1,1){\sgl}}
%\put(\sgl,\sgl){\line(-1,0){\sgl}}
%\put(\sgl,\sgl){\line(0,-1){\sgl}}
\put(0,\sgl){\line(1,-1){\sgl}}
\put(-8,0){1}
\put(23,0){2}
\put(-8,20){3}
%\put(23,\gl){4}
\put(0,0){\circle*{2}}
\put(\sgl,0){\circle*{4}}
\put(0,\sgl){\circle*{2}}
%\put(\sgl,\sgl){\circle*{2}}
%\put(\neghalfsgl,\halfsgl){24}
\put(0,\negsgl){32}
\end{picture}}
\end{picture}
\caption{Labeled rooted trees on three vertices: Pr\"ufer sequences}
\end{figure}
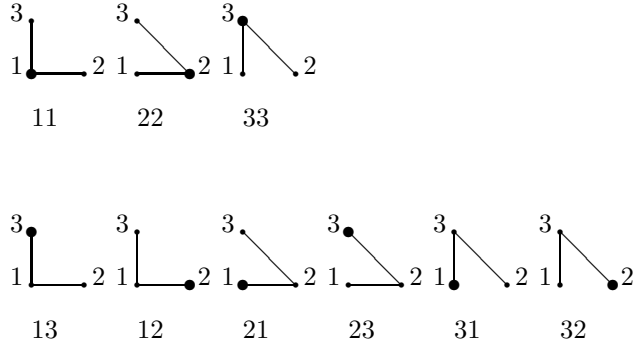

%Every forest on a set $V$ with $m$ elements comes from a   tree on a label set $U$ with one more point. Since there are
%$n^{n-1}$ labeled rooted trees, there are $n^{n-2}$ labeled  trees. Hence there are $(m+1)^{m-1}$ forests.

\noindent\textbf{Labeled rooted trees defined recursively}
A labeled rooted tree $T$ on $U$ has a recursive definition as a point $r$ in $U$ together with a forest of labeled rooted trees on $U \setminus \{r\}$.
This forest $F$ is the restriction of $T$ to points in $U \setminus \{r\}$ that are not in $T^{-1}(r)$. The set of roots of the forest is $R = T^{-1}(r)$.
The recursive definition ends when the tree consists only of a root; the forest is then empty.

\noindent\textbf{Notation for labeled rooted trees}
A labeled rooted tree on a one-point set may be designated by its label $j$. A labeled rooted tree on a set with two or more points may be
denoted by $j [ \; \; \; \;]$, where the forest of immediate successor rooted trees is listed (in some arbitrary order) inside the bracket. As an example, consider the
tree with root $c$ and with vertices $b,e$ that are sent to $c$ and vertices $a,d$ that are sent to $e$. In this notation the tree
would be $c[be[ad]]$.

\noindent\textbf{Labeled rooted trees as sequences of vertices} Consider a non-empty vertex set $U$ with $n$ elements. For the construction to follow it
is necessary to impose a linear order on $U$. The result says that labeled rooted trees on $U$ correspond to sequences of $n-1$ elements.

\begin{example} For $n=3$ the vertices of a labeled rooted tree may be numbered 1,2,3. As shown below, each labeled rooted tree may be
 coded by a sequence of two numbers.
For example, the tree 2[13] is coded by 22, wile the tree $2[1[3]]$ is coded by 12. There are nine sequences: three of them 11, 22, 33 correspond to one unlableled rooted  tree, while six of them 12, 13, 21, 23, 31, 33 correspond to the other unlabeled rooted tree.
See Figure~2 for the picture.
\end{example}

\begin{proposition}[Pr\"ufer correspondence] Given non-empty label set $U$ with a given linear order,
there is a bijection between the set of
sequences $s$ of length $n-1$ of elements of $U$ and the set of labeled rooted trees $T$ in $\Ac[U]$.
For each $j$ in $U$ the number of times the sequence $s$ assumes the value $j$ is
the degree $|T^{-1}(j)|$.
\end{proposition}

\begin{proof}It is easiest to see how to go from the labeled rooted tree $T$ to the corresponding sequence $s$. At each stage
remove the smallest leaf and its corresponding edge from the tree. The value of $s$ at this
stage is the vertex at the other end of this edge. When this is repeated $k-1$ times the final
sequence value is the root. This procedure is illustrated in Figure~3.

Here is the construction to go from the sequence $s$ to the labeled rooted tree. The edges are restored in the same
order as they were removed. At each stage add a new edge as follows.
Take the smallest vertex that has not yet been used and that does not occur in the part of the sequence that has not
been used. The edge then goes from this vertex to the next element of the sequence. See Figure~4 for a picture.
\end{proof}

\begin{proposition}[Cayley] The number of labeled rooted trees $T$ with
$|T| = n$ vertices is $n^{n-1}$.
\end{proposition}

\begin{proof}
This famous result of Cayley results from the Pr\"ufer correspondence. There are $n^{n-1}$ sequences of length $n-1$ in a set $U$ with $n$ elements.
\end{proof}

% Figure 3 From rooted tree to Prufer sequence
\begin{figure}[tb]
\begin{picture}(\bgsgl,\foursgl)(\negtwosgl,\sgl)
\put(0,\twosgl){
\begin{picture}(\sgl,\sgl)
\put(0,0){\line(1,0){\sgl}}
\put(0,0){\line(0,1){\sgl}}
\put(0,\sgl){\line(1,0){\sgl}}
\put(0,\sgl){\line(0,1){\sgl}}
\put(0,0){\circle*{5}}
\put(0,\sgl){\circle*{2}}
\put(0,\twosgl){\circle*{2}}
\put(\sgl,0){\circle*{2}}
\put(\sgl,\sgl){\circle*{2}}
\put(\neghalfsgl,0){3}
\put(\neghalfsgl,\sgl){1}
\put(\neghalfsgl,\twosgl){2}
\put(\sgl,0){$\;4$}
\put(\sgl,\sgl){$\;5$}
%\put(0,\negsgl){1313}
\end{picture}}
\put(\threesgl,\twosgl){
\begin{picture}(\sgl,\sgl)
\put(0,0){\line(1,0){\sgl}}
\put(0,0){\line(0,1){\sgl}}
\put(0,\sgl){\line(1,0){\sgl}}
%\put(0,\sgl){\line(0,1){\sgl}}
\put(0,0){\circle*{5}}
\put(0,\sgl){\circle*{2}}
\put(0,\twosgl){\circle*{2}}
\put(\sgl,0){\circle*{2}}
\put(\sgl,\sgl){\circle*{2}}
\put(\neghalfsgl,0){3}
\put(\neghalfsgl,\sgl){1}
\put(\neghalfsgl,\twosgl){2}
\put(\sgl,0){$\;4$}
\put(\sgl,\sgl){$\;5$}
\put(0,\negsgl){1}
\end{picture}}
\put(\sixsgl,\twosgl){
\begin{picture}(\sgl,\sgl)
%\put(0,0){\line(1,0){\sgl}}
\put(0,0){\line(0,1){\sgl}}
\put(0,\sgl){\line(1,0){\sgl}}
%\put(0,\sgl){\line(0,1){\sgl}}
\put(0,0){\circle*{5}}
\put(0,\sgl){\circle*{2}}
\put(0,\twosgl){\circle*{2}}
\put(\sgl,0){\circle*{2}}
\put(\sgl,\sgl){\circle*{2}}
\put(\neghalfsgl,0){3}
\put(\neghalfsgl,\sgl){1}
\put(\neghalfsgl,\twosgl){2}
\put(\sgl,0){$\;4$}
\put(\sgl,\sgl){$\;5$}
\put(0,\negsgl){13}
\end{picture}}
\put(\ninesgl,\twosgl){
\begin{picture}(\sgl,\sgl)
%\put(0,0){\line(1,0){\sgl}}
\put(0,0){\line(0,1){\sgl}}
%\put(0,\sgl){\line(1,0){\sgl}}
%\put(0,\sgl){\line(0,1){\sgl}}
\put(0,0){\circle*{5}}
\put(0,\sgl){\circle*{2}}
\put(0,\twosgl){\circle*{2}}
\put(\sgl,0){\circle*{2}}
\put(\sgl,\sgl){\circle*{2}}
\put(\neghalfsgl,0){3}
\put(\neghalfsgl,\sgl){1}
\put(\neghalfsgl,\twosgl){2}
\put(\sgl,0){$\;4$}
\put(\sgl,\sgl){$\;5$}
\put(0,\negsgl){131}
\end{picture}}
\put(\twelvesgl,\twosgl){
\begin{picture}(\sgl,\sgl)
%\put(0,0){\line(1,0){\sgl}}
%\put(0,0){\line(0,1){\sgl}}
%\put(0,\sgl){\line(1,0){\sgl}}
%\put(0,\sgl){\line(0,1){\sgl}}
\put(0,0){\circle*{5}}
\put(0,\sgl){\circle*{2}}
\put(0,\twosgl){\circle*{2}}
\put(\sgl,0){\circle*{2}}
\put(\sgl,\sgl){\circle*{2}}
\put(\neghalfsgl,0){3}
\put(\neghalfsgl,\sgl){1}
\put(\neghalfsgl,\twosgl){2}
\put(\sgl,0){$\;4$}
\put(\sgl,\sgl){$\;5$}
\put(0,\negsgl){1313}
\end{picture}}
\end{picture}
\caption{From labeled rooted tree to Pr\"ufer sequence}
\end{figure}
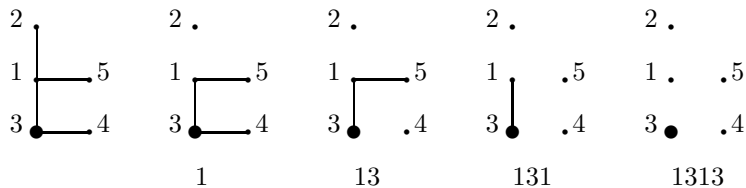

% Figure 4 From Prufer sequence to rooted tree
\begin{figure}[tb]
\begin{picture}(\bgsgl,\threesgl)(\negtwosgl,\sgl)
\put(0,\twosgl){
\begin{picture}(\sgl,\sgl)
%\put(0,0){\line(1,0){\sgl}}
%\put(0,0){\line(0,1){\sgl}}
%\put(0,\sgl){\line(1,0){\sgl}}
%\put(0,\sgl){\line(0,1){\sgl}}
\put(0,0){\circle*{2}}
\put(0,\sgl){\circle*{2}}
\put(0,\twosgl){\circle*{2}}
\put(\sgl,0){\circle*{2}}
\put(\sgl,\sgl){\circle*{2}}
\put(\neghalfsgl,0){3}
\put(\neghalfsgl,\sgl){1}
\put(\neghalfsgl,\twosgl){2}
\put(\sgl,0){$\;4$}
\put(\sgl,\sgl){$\;5$}
\put(0,\negsgl){1313}
\end{picture}}
\put(\threesgl,\twosgl){
\begin{picture}(\sgl,\sgl)
%\put(0,0){\line(1,0){\sgl}}
%\put(0,0){\line(0,1){\sgl}}
%\put(0,\sgl){\line(1,0){\sgl}}
\put(0,\sgl){\line(0,1){\sgl}}
\put(0,0){\circle*{2}}
\put(0,\sgl){\circle*{2}}
\put(0,\twosgl){\circle*{2}}
\put(\sgl,0){\circle*{2}}
\put(\sgl,\sgl){\circle*{2}}
\put(\neghalfsgl,0){3}
\put(\neghalfsgl,\sgl){1}
\put(\neghalfsgl,\twosgl){2}
\put(\sgl,0){$\;4$}
\put(\sgl,\sgl){$\;5$}
\put(0,\negsgl){\ 313}
\end{picture}}
\put(\sixsgl,\twosgl){
\begin{picture}(\sgl,\sgl)
\put(0,0){\line(1,0){\sgl}}
%\put(0,0){\line(0,1){\sgl}}
%\put(0,\sgl){\line(1,0){\sgl}}
\put(0,\sgl){\line(0,1){\sgl}}
\put(0,0){\circle*{2}}
\put(0,\sgl){\circle*{2}}
\put(0,\twosgl){\circle*{2}}
\put(\sgl,0){\circle*{2}}
\put(\sgl,\sgl){\circle*{2}}
\put(\neghalfsgl,0){3}
\put(\neghalfsgl,\sgl){1}
\put(\neghalfsgl,\twosgl){2}
\put(\sgl,0){$\;4$}
\put(\sgl,\sgl){$\;5$}
\put(0,\negsgl){\ \ 13}
\end{picture}}
\put(\ninesgl,\twosgl){
\begin{picture}(\sgl,\sgl)
\put(0,0){\line(1,0){\sgl}}
%\put(0,0){\line(0,1){\sgl}}
\put(0,\sgl){\line(1,0){\sgl}}
\put(0,\sgl){\line(0,1){\sgl}}
\put(0,0){\circle*{2}}
\put(0,\sgl){\circle*{2}}
\put(0,\twosgl){\circle*{2}}
\put(\sgl,0){\circle*{2}}
\put(\sgl,\sgl){\circle*{2}}
\put(\neghalfsgl,0){3}
\put(\neghalfsgl,\sgl){1}
\put(\neghalfsgl,\twosgl){2}
\put(\sgl,0){$\;4$}
\put(\sgl,\sgl){$\;5$}
\put(0,\negsgl){\ \ \ 3}
\end{picture}}
\put(\twelvesgl,\twosgl){
\begin{picture}(\sgl,\sgl)
\put(0,0){\line(1,0){\sgl}}
\put(0,0){\line(0,1){\sgl}}
\put(0,\sgl){\line(1,0){\sgl}}
\put(0,\sgl){\line(0,1){\sgl}}
\put(0,0){\circle*{5}}
\put(0,\sgl){\circle*{2}}
\put(0,\twosgl){\circle*{2}}
\put(\sgl,0){\circle*{2}}
\put(\sgl,\sgl){\circle*{2}}
\put(\neghalfsgl,0){3}
\put(\neghalfsgl,\sgl){1}
\put(\neghalfsgl,\twosgl){2}
\put(\sgl,0){$\;4$}
\put(\sgl,\sgl){$\;5$}
%\put(0,\negsgl){1313}
\end{picture}}
\end{picture}
\caption{From Pr\"ufer sequence to labeled rooted  tree}
\end{figure}
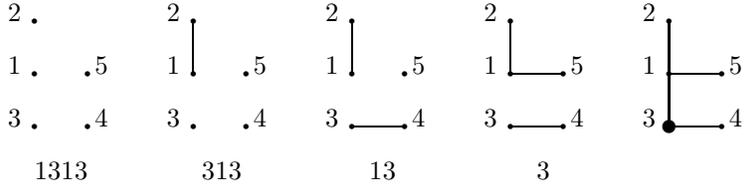

\section{Unlabeled rooted trees}

The difficulty with unlabeled rooted trees is the general difficulty with unlabeled combinatorial structures (isomorphism classes
of structures). This is the presence of symmetry. This is a well-studied topic; there are nice accounts in \cite{BergeronLabelleLeroux} and
in \cite{Kerber}.

\noindent\textbf{Unlabeled rooted trees via orbits of labeled rooted trees} An \emph{unlabeled rooted tree} is an isomorphism
invariant of labeled rooted trees. It is convenient to fix the label set $U$.
For each labeled rooted tree $T$ on $U$ there is a corresponding unlabeled rooted tree $\tau$.  Here are the details.

Let $U$ be a non-empty set. Let $\Ac[U]$ be the set of labeled rooted trees on vertex set $U$. Each such tree is a
function $T: U \setminus \{r\}$ to $U$. For each $T$ in
$\Ac[U]$ and each bijection $b: U \to U$, the composite function $T \circ b$ is   another labeled rooted tree in $\Ac[U]$.
It is a function $T' = T \circ b: U' \setminus \{r'\}$ to $U$, where $b r' = r$. Two such labeled rooted trees $T, T'$ are
said to be isomorphic.

An \emph{unlabeled rooted tree} $\tau$ with $n$ vertices is an object that corresponds to an isomorphism class of labeled rooted trees,
where the label set has $n$ elements. The set of unlabeled rooted trees with $n$ vertices is denoted $\tilde \Ac[n]$. The set
of all unlabeled rooted trees is denoted $\tilde \Ac$. There is no unlabeled rooted tree with zero vertices, but sometimes it is
convenient to introduce an extra empty set object associated with zero vertices. The augmented set is written $\tilde \Ac_\emptyset$.

An unlabeled rooted tree $\tau$ has no underlying set and thus no vertices and no edges. Nevertheless,
it may be pictured by any $T$ in the isomorphism class.
There are various invariants of an unlabeled rooted tree
$\tau$. Among them are
the number of vertices $|\tau|$ and the number of vertices $v_k(\tau)$ of degree $k$. These are related by
\begin{equation}
\sum_k k v_k(\tau) = |\tau| - 1.
\end{equation}

There are relatively few unlabeled rooted trees. The simplest (but nevertheless important) ones have 1 root and $n-1$ leaves. In
the following such a tree will be denoted $\tau = n-1$. For this tree $|\tau| = n$ and $v_0(\tau) = n-1$, $v_{n-1}(\tau) = 1$.
Another simple class of unlabeled rooted trees are the linear trees with $v_0(\tau) = 1$ and $v_1(\tau) = n-1$.

Group theory illuminates the situation. Fix a label set $U$
with $n$ elements. Let $\Gc$ be the permutation group of this set. This consists of all bijections $b: U \to U$.  This group has $|\Gc| = n!$ elements.
The group $\Gc$ also acts on the set $\Ac[U]$ of labeled rooted trees with  $n^{n-1}$ elements. For each bijection $b$, the map $T \to T \circ b$ is a map from $\Ac[U]$ to itself.
The set $\tilde \Ac[n]$ of unlabeled rooted trees corresponds to the set of orbits under this action:
\begin{equation}
\tilde \Ac[n] \cong \Ac[U]/\Gc.
\end{equation}
For each labeled rooted tree $T$, the corresponding unlabeled rooted
tree  $\tau$ is an abstract object corresponding to the orbit $\Gc T$ of $T$.
The map from labeled rooted trees to unlabeled rooted trees may be
summarized by the surjection
\begin{equation}
\Ac[U] \to \tilde\Ac[n].
\end{equation}

According to the theory of group actions,  the size of the orbit $\Gc T$ of $T$ is the order $|\Gc| = n!$   divided by
the order $|\Gc_T|$, where $\Gc_T$ is the stabilizer subgroup of $T$. Thus
\begin{equation}
|\Gc T| = \frac{n!}{|\Gc_T|} .
\end{equation}
The order $|\Gc_T|$ is the same for all $T$ in an orbit and hence
may be denoted $\sigma(\tau)$, where $\tau$ is the unlabeled rooted tree corresponding to $T$. The number $\sigma(\tau)$ is
the \emph{symmetry factor} of $\tau$. The size of the orbit also depends only on $\tau$; it will be denoted $r(\tau)$.
This gives the following basic result:

\begin{proposition} For rooted trees with $|\tau| = n$ vertices the number of labeled rooted trees in $\Ac[U]$ per
unlabeled rooted tree in $\tilde \Ac[n]$ is
\begin{equation}
r(\tau) = \frac{n!}{\sigma(\tau)} .
\end{equation}
\end{proposition}

\begin{example} There are labeled rooted trees on $n=4$ vertices that consist of a root and three leaves. The orbit $\Gc T$ of such a rooted tree is shown in Figure~5. For each rooted tree $T$ in this orbit, the corresponding stabilizer subgroup $\Gc_T$ has 6 elements, corresponding to the 6 permutations of the leaves. The symmetry factor is $\sigma(\tau) = 6$. The number of labeled rooted trees in the orbit is $r(\tau) =  24/6 = 4$.
\end{example}

% Figure 5 Orbit of a labeled rooted tree on 4 vertices
\begin{figure}[tb]
\begin{picture}(\bgsgl,\twosgl)(\negfoursgl,0)
\put(0,0){
\begin{picture}(\sgl,\sgl)
\put(0,0){\line(1,0){\sgl}}
\put(0,0){\line(0,1){\sgl}}
\put(0,0){\line(1,1){\sgl}}
%\put(\sgl,\sgl){\line(-1,0){\sgl}}
%\put(\sgl,\sgl){\line(0,-1){\sgl}}
%\put(0,\sgl){\line(1,-1){\sgl}}
\put(0,0){\circle*{4}}
\put(\sgl,0){\circle*{2}}
\put(0,\sgl){\circle*{2}}
\put(\sgl,\sgl){\circle*{2}}
%\put(\neghalfsgl,\halfsgl){4}
\end{picture}}
\put(\twosgl,0){
\begin{picture}(0,0)
\put(0,0){\line(1,0){\sgl}}
%\put(0,0){\line(0,1){\sgl}}
%\put(0,0){\line(1,1){\sgl}}
%\put(\sgl,\sgl){\line(-1,0){\sgl}}
\put(\sgl,\sgl){\line(0,-1){\sgl}}
\put(0,\sgl){\line(1,-1){\sgl}}
\put(0,0){\circle*{2}}
\put(\sgl,0){\circle*{4}}
\put(0,\sgl){\circle*{2}}
\put(\sgl,\sgl){\circle*{2}}
%\put(\neghalfsgl,\halfsgl){24}
\end{picture}}
\put(\foursgl,0){
\begin{picture}(\sgl,\sgl)
%\put(0,0){\line(1,0){\sgl}}
%\put(0,0){\line(0,1){\sgl}}
\put(0,0){\line(1,1){\sgl}}
\put(\sgl,\sgl){\line(-1,0){\sgl}}
\put(\sgl,\sgl){\line(0,-1){\sgl}}
%\put(0,\sgl){\line(1,-1){\sgl}}
\put(0,0){\circle*{2}}
\put(0,\sgl){\circle*{2}}
\put(\sgl,0){\circle*{2}}
\put(\sgl,\sgl){\circle*{4}}
%\put(\sgl,\sgl){\circle*{2}}
%\put(\neghalfsgl,\halfsgl){12}
\end{picture}}
\put(\sixsgl,0){
\begin{picture}(\sgl,\sgl)
%\put(0,0){\line(1,0){\sgl}}
\put(0,0){\line(0,1){\sgl}}
%\put(0,0){\line(1,1){\sgl}}
\put(\sgl,\sgl){\line(-1,0){\sgl}}
%\put(\sgl,\sgl){\line(0,-1){\sgl}}
\put(0,\sgl){\line(1,-1){\sgl}}
\put(0,0){\circle*{2}}
\put(\sgl,0){\circle*{2}}
\put(0,\sgl){\circle*{4}}
\put(\sgl,\sgl){\circle*{2}}
%\put(\neghalfsgl,\halfsgl){24}
\end{picture}}
\end{picture}
\caption{Orbit of a labeled rooted tree on four vertices: 4 trees, $\sigma(\tau) = 6$}
\end{figure}
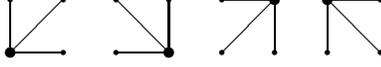

There is an identity that expresses the fact that the sum over unlabeled rooted trees of
the corresponding number of labeled rooted trees is the total number of labeled rooted trees.

\begin{proposition}
For unlabeled rooted trees $\tau$ with $|\tau| = n$ vertices the sum of the corresponding numbers $r(\tau)$ of labeled rooted trees in
$\Ac[U]$ with $|U| = n$ is
\begin{equation}
\sum_{\tau \in \tilde\Ac[n]} \frac{n!}{\sigma(\tau)} = n^{n-1}.
\end{equation}
\end{proposition}

\noindent\textbf{Multisets of unlabeled rooted trees} The analog of a forest of labeled rooted trees is a multiset of unlabeled rooted trees.
A multiset of unlabeled rooted trees with $m$ vertices is defined as a function $N$ from unlabeled rooted trees in $\tilde \Ac$ to natural numbers $\geq 0$
such that
\begin{equation}
\sum_{\tau} |\tau| N(\tau) = m.
\end{equation}
 Thus $N(\tau)$ represents the number of times the unlabeled rooted tree $\tau$ occurs in the multiset.
 Every forest of labeled rooted trees gives rise to a multiset of unlabeled rooted trees, where $N(\tau)$ is the number of blocks
 in the forest that correspond to unlabeled rooted tree $\tau$.

 Such a multiset has a \emph{symmetry factor} derived from the symmetry factor associated with  unlabeled rooted trees.
  Let $F$ be
 a forest with corresponding multiset $N$. The symmetry factor $\sigma(N)$ is
  the order of the stabilizer subgroup $\Gc_F$ of the forest. It is
 \begin{equation}
 \sigma(N) = \prod_{\tau'}  N(\tau')!\sigma(\tau')^{N(\tau')} .
 \end{equation}
 This expression may be derived using group theory.
 Let $\Bc$ be the subgroup of $\Gc_F$ generated by permutations
 that leave each block invariant. For each block there is a corresponding symmetry factor $\sigma(\tau')$, so $\Bc$ is  a product group with order $\prod_{\tau'} \sigma(\tau')^{N(\tau')}$.
 Let $\Hc$ be the subgroup that permutes blocks with identical unlabeled rooted trees. If unlabeled rooted tree $\tau'$ occurs in $N(\tau')$ blocks, then
 there are $N(\tau')!$ permutations involving that rooted tree. The order of $\Hc$ is thus $\prod_{\tau'} N(\tau')! $. The group $\Hc$ is a normal subgroup of $\Gc_F$,
 so in particular for every $b'$ in $\Bc$ and for every $h$ in $\Hc$ the element $h b' h^{-1}$ is also in $\Bc$.
 Every element of $\Gc_F$ may be uniquely expressed as  product $bh$ of an element of $\Bc$ with an element of $\Hc$. (This decomposition respects
 multiplication: $bh b'h' = (b h b' h^{-1}) (h h')$. In fact the group $\Gc_F$ is the
 semidirect product of the group $\Bc$ with the group $\Hc$ \cite{MacLaneBirkhoff}.)  The conclusion is that the order of $\Gc_F$ is the product of the orders of $\Bc$ and $\Hc$.

Here is how to construct a forest corresponding to a given multiset $N$ with $m$ vertices.
Find a set $V$ with $m$ elements and a set partition $\Gamma$ of $V$. Require that there is a function
$\chi: \Gamma \to \Ac$ such that for each $\tau$ the inverse image $\chi^{-1}(\tau)$
consists of $N(\tau)$ blocks of size $|\tau|$. Finally, for each block $B$ in $\Gamma$  find a labeled rooted tree $T$ that
determines unlabeled rooted tree $\chi(B)$.

The number of pairs $\Gamma, \chi$ satisfying these conditions is the coefficient
 \begin{equation}
 C(N) = \frac{m!}{\prod_{\tau'} (|\tau'|!)^{N(\tau')} } \frac{1}{\prod_{\tau'} N(\tau')!}.
 \end{equation}
The first factor is the multinomial coefficient that determines how many ways of producing blocks
of the appropriate sizes in a given order. The second factor has a denominator that describes how many
ways there are of permuting the blocks to preserve $\chi$. The number of forests
is
 \begin{equation}
 f(N) = C(N)  \prod_{\tau'}  r(\tau')^{N(\tau')} = \frac{m!}{\prod_{\tau'} N(\tau')!\sigma(\tau')^{N(\tau')} }=
 \frac{m!}{\sigma(N)}.
 \end{equation}
The first equality comes from counting  the number of
ways of putting labeled rooted trees in the appropriate blocks. The second equality results from inserting $r(\tau') = |\tau'|!/\sigma(\tau')$.

\noindent\textbf{Unlabeled rooted trees defined recursively}
There is a recursive definition of unlabeled rooted trees.
An unlabeled rooted tree $\tau$ with $n\geq 1$ vertices is equivalent to a multiset $N$ of unlabeled rooted trees
with  $n-1$ vertices. This counts the  unlabeled rooted subtrees that result when the root is removed.
The recursion terminates with unlabeled rooted trees with one vertex; the corresponding multiset is zero.

The number $r(\tau)$ that counts labeled rooted trees satisfies the recursion
\begin{equation}
r(\tau) =  |\tau|  f(N)  = |\tau| C(N) \prod_{\tau'}  r(\tau')^{N(\tau')} .
\end{equation}
Here $\tau$ is an unlabeled rooted tree on $n$ vertices, and  $N$ is the corresponding multiset on $n-1$ vertices. This is
because  a labeled rooted tree is determined by a root point and a forest over the remaining points.

There is also a recursion relation for the symmetry factors. If $\tau$ has $n$ vertices and its subtrees define a multiset $N$ with $n-1$ vertices.
\begin{equation}
\sigma(\tau) = \sigma(N) =  \prod_{\tau' \in \tilde \Ac}  N(\tau')!\sigma(\tau')^{N(\tau')}.
\end{equation}

This recursion relation has an explicit solution. Let $\tau$ be an unlabeled rooted tree. Consider some labeling, so
that there is a set $U$ of vertices and a labeled rooted tree on $U$. For each vertex $j$, consider the subtree above $j$, and let $N_j$ count the
unlabeled rooted trees above this subtree. Then
\begin{equation}
\sigma(\tau) = \prod_j \prod_{\tau'} N_j(\tau')!
\end{equation}

\noindent\textbf{Notation for unlabeled rooted trees}
The multiset notation gives a convenient way of describing unlabeled rooted trees. A tree with a single vertex is denoted 0.
Otherwise, the tree is denoted $N_1[\tau_1] N_2[\tau_2] \ldots N_k[\tau_k]$, where each $N_j \neq 0$ and the $\tau_j$ are
descriptions of different unlabeled rooted trees. It is convenient to abbreviate $N[0]$ by $N$. Thus, for example, the
labeled rooted tree $c[b[e[ad]]$ would determine the unlabeled rooted tree $1[0]1[2[0]]$. In the abbreviated form this would be
$11[2]$. This says that the root has 1 immediate successor with a single vertex and 1 immediate successor that is a tree with 2
immediate successor vertices.

\begin{example} For $n=1$ the only rooted tree is 0, consisting of a single root point.
For  $n=2$ and a given vertex set there are two labeled rooted trees, depending on which point is chosen for the root.
There is only one unlabeled rooted tree, denoted 1.
For $n=3$ and a given vertex set there
 are $3^2 = 9$ labeled rooted trees. These decompose into two orbits, as shown in Figure~1. These
 correspond to unlabeled rooted trees $\tau$ that may be denoted $2$ and $1[1]$. The rooted tree 2 has a root with two  leaves. The symmetry factor is 2.
The $1[1]$ linear rooted tree has a root and a successor rooted tree $1$. The symmetry factor is 1.
This gives the correct number of labeled rooted trees as the sum
$6/2 + 6/1 = 9 = 3^2$. The two unlabeled rooted trees are shown in Figure~6.
\end{example}

% Figure 6 Unlabeled rooted trees on 3 vertices
\begin{figure}[tb]
\begin{picture}(\bggl,\gl)(\negthreegl,0)
\put(0,0){
\begin{picture}(\gl,\gl)
\put(0,0){\line(1,0){\gl}}
\put(0,0){\line(0,1){\gl}}
%\put(0,0){\line(1,1){\gl}}
%\put(\gl,\gl){\line(-1,0){\gl}}
%\put(\gl,\gl){\line(0,-1){\gl}}
%\put(0,\gl){\line(1,-1){\gl}}
\put(0,0){\circle*{2}}
\put(0,\gl){\circle*{2}}
\put(\gl,0){\circle*{2}}
%\put(\gl,\gl){\circle*{2}}
\put(0,0){\circle*{4}}
\put(\neghalfgl,\halfgl){3}
\end{picture}}
\put(\twogl,0){
\begin{picture}(0,0)
\put(0,0){\line(1,0){\gl}}
%\put(0,0){\line(0,1){\gl}}
%\put(0,0){\line(1,1){\gl}}
%\put(\gl,\gl){\line(-1,0){\gl}}
%\put(\gl,\gl){\line(0,-1){\gl}}
\put(0,\gl){\line(1,-1){\gl}}
\put(0,0){\circle*{2}}
\put(0,\gl){\circle*{2}}
\put(\gl,0){\circle*{2}}
%\put(\gl,\gl){\circle*{2}}
\put(0,0){\circle*{4}}
\put(\neghalfgl,\halfgl){6}
\end{picture}}
\end{picture}
\caption{Unlabeled rooted trees on three vertices: $\sigma(\tau) = 2,1$}
\end{figure}
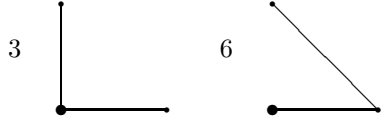

\begin{example}
The case $n=4$ is more interesting. There are  four unlabeled rooted trees, which may be denoted in multiset notation by
3, 11[1], 1[2], and 1[1[1]].  These rooted trees
 have   symmetry factors $6,1,2,1$.
The number of labeled rooted trees is the sum
$24/6 + 24/1 + 24/2 + 24/1 = 64 = 4^3$. See Figure~7 for a picture of the unlabeled rooted trees.
\end{example}

% Figure 7 Unlabeled rooted trees on 4 vertices
\begin{figure}[tb]
\begin{picture}(\bggl,\gl)(\neggl,0)
\put(0,0){
\begin{picture}(\gl,\gl)
\put(0,0){\line(1,0){\gl}}
\put(0,0){\line(0,1){\gl}}
\put(0,0){\line(1,1){\gl}}
%\put(\gl,\gl){\line(-1,0){\gl}}
%\put(\gl,\gl){\line(0,-1){\gl}}
%\put(0,\gl){\line(1,-1){\gl}}
\put(0,0){\circle*{2}}
\put(0,\gl){\circle*{2}}
\put(\gl,0){\circle*{2}}
\put(\gl,\gl){\circle*{2}}
\put(0,0){\circle*{4}}
\put(\neghalfgl,\halfgl){4}
\end{picture}}
\put(\twogl,0){
\begin{picture}(0,0)
\put(0,0){\line(1,0){\gl}}
\put(0,0){\line(0,1){\gl}}
%\put(0,0){\line(1,1){\gl}}
%\put(\gl,\gl){\line(-1,0){\gl}}
\put(\gl,\gl){\line(0,-1){\gl}}
%\put(0,\gl){\line(1,-1){\gl}}
\put(0,0){\circle*{2}}
\put(0,\gl){\circle*{2}}
\put(\gl,0){\circle*{2}}
\put(\gl,\gl){\circle*{2}}
\put(0,0){\circle*{4}}
\put(\neghalfgl,\halfgl){24}
\end{picture}}
\put(\fourgl,0){
\begin{picture}(\gl,\gl)
\put(0,0){\line(1,0){\gl}}
%\put(0,0){\line(0,1){\gl}}
%\put(0,0){\line(1,1){\gl}}
%\put(\gl,\gl){\line(-1,0){\gl}}
\put(\gl,\gl){\line(0,-1){\gl}}
\put(0,\gl){\line(1,-1){\gl}}
\put(0,0){\circle*{2}}
\put(0,\gl){\circle*{2}}
\put(\gl,0){\circle*{2}}
\put(\gl,\gl){\circle*{2}}
\put(0,0){\circle*{4}}
\put(\neghalfgl,\halfgl){12}
\end{picture}}
\put(\sixgl,0){
\begin{picture}(\gl,\gl)
\put(0,0){\line(1,0){\gl}}
%\put(0,0){\line(0,1){\gl}}
%\put(0,0){\line(1,1){\gl}}
\put(\gl,\gl){\line(-1,0){\gl}}
\put(\gl,\gl){\line(0,-1){\gl}}
%\put(0,\gl){\line(1,-1){\gl}}
\put(0,0){\circle*{2}}
\put(0,\gl){\circle*{2}}
\put(\gl,0){\circle*{2}}
\put(\gl,\gl){\circle*{2}}
\put(0,0){\circle*{4}}
\put(\neghalfgl,\halfgl){24}
\end{picture}}
\end{picture}
\caption{Unlabeled rooted trees on four vertices: $\sigma(\tau) = 6, 1, 2, 1$}
\end{figure}
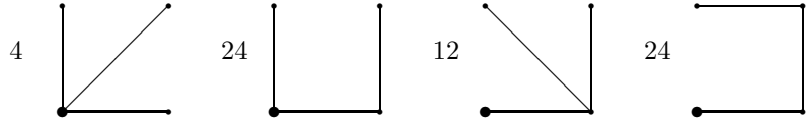

\begin{remark}
The formula for the number $a_n$ of labeled rooted trees on a vertex set with $n$ vertices is $a_n = n^{n-1}$. The
 number $\tilde a_n$ of unlabeled rooted trees with $n$ vertices is not so easy to compute. This may be
seen by contrasting the generating functions.

The exponential generating function for the number $a_n$ of labeled rooted trees with $n$ vertices is
$a(t) = \sum_{n = 1}^\infty \frac{1}{n!} a_n t^n$.  The recursive definition of labeled rooted tree
from root point and forest gives
\begin{equation}
a(t) = t \exp( a(t)) .
\end{equation}
Labeled enumeration give a simple result:  for fixed $t$ the value $x= a(t)$ satisfies the fixed point equation $x = t \exp(x)$.

Contrast this with the  generating function for the number $\tilde a_n$ of unlabeled rooted trees with $n$ vertices. This is
$\tilde a(t) = \sum_{n = 1}^\infty \tilde a_n t^n$.
 The recursive definition of unlabeled rooted tree from multiset
gives the identity
\begin{equation}
\tilde a(t) = t  \exp( \sum_{k=1}^\infty \frac{1}{k} \tilde a(t^k)) .
\end{equation}
Unlabeled enumeration produces a much more complicated equation.
See \cite{FlajoletSedgewick} or \cite{BergeronLabelleLeroux} for the full story.
\end{remark}

\section{Fixed point equations and labeled rooted trees}

Let $\beta(x)$ be a formal power series in $x$. Let $t$ and $g$ be parameters. Consider
the \emph{fixed point equation}
\begin{equation}
x = g + t \beta(x).
\end{equation}

\begin{proposition} The fixed point equation $x = g + t \beta(x)$ has the formal solution
\begin{equation}
x = f(t,g) = \sum_{n=0}^\infty \frac{t^n}{n!} f_n(g) ,
\end{equation}
where $f_0(g) = g$, and where for for $n \geq 1$ the coefficient $f_n(g)$ has the explicit representation
\begin{equation}
f_n(g) = \left( \frac{\partial}{\partial g} \right)^{n-1} \beta(g)^n.
\end{equation}
\end{proposition}

\begin{proof}
Fix $g$. The problem is to find the expansion of $x$ as a function of $t$. We know the inverse function
\begin{equation}
t = (x-g)/\beta(x)
\end{equation} giving $t$ as a function of $x$. Notice that $t=0$ corresponds to $x=g$.
The Lagrange inversion formula applies. The formula is based on the
the fact that the residue of a differential form expressed by a formal Laurent series is invariant under change of variable. Start with the
identity
\begin{equation}
\frac{1}{n} d \left(\frac{x}{t^n} \right) =  - \frac{x}{t^{n+1}}  \, dt + \frac{1}{n} \frac{1}{t^n}  \, dx.
\end{equation}
Since the left hand is a perfect differential, it has residue zero. This gives an identity for the residues
\begin{equation}
\residue \left(\frac{x}{t^{n+1}}  \, dt\right) = \frac{1}{n} \residue \left(\frac{1}{t^n}  \, dx \right).
\end{equation}
This is the Lagrange inversion formula.

In the case at hand
\begin{equation}
\frac{1}{n!} f_n(g) = \residue \left(\frac{x}{t^{n+1}}  \, dt\right) = \frac{1}{n} \residue \left(\frac{1}{t^n}  \, dx \right) =
\frac{1}{n} \residue \left(\frac{\beta(x)^n}{(x-g)^n}  \, dx \right),
\end{equation}
where the last residue is computed at the singularity $x=g$.
The residue is $1/(n-1)!$ times the $n-1$st derivative of $\beta(g)^n$ with respect to $g$. This gives the result.
\end{proof}

A combinatorial solution gives more detailed information.
This is given by an expansion indexed by rooted trees.
For each $n\geq 1$  fix a label set $U_n$ with $n$ elements, and  consider rooted trees $T$ with the label set as vertex set.
For $n = 0$ there is an empty set object associated with the label set $U_0 = \emptyset$.

\begin{proposition} The fixed point equation $x = g + t \beta(x)$ has the  solution
given by a formal power series as above,
where for $n \geq 1$ the coefficient $f_n(g)$ has the explicit representation
\begin{equation}
f_n(g) = \sum_{T \in \Ac[U_n]} f_T(g)
\end{equation}
For $n \geq 1$ the coefficient
\begin{equation}
f_T(g) = \prod_{j \in U_n} \left( \frac{\partial}{\partial g} \right)^{|T^{-1}(j)|} \beta(g),
\end{equation}
where $|T^{-1}(j)|$ is the degree of vertex $j$ of rooted tree $T$. For $n=0$ the contribution from
the empty rooted tree is $f_\emptyset(g) = g$.
\end{proposition}

\begin{proof} This proof uses calculus formula in the form explained in the appendix.
Identify the $n$ factors in $\beta(g)^n$ with the $n$ points in the vertex set $U_n$. Each derivative corresponds to
a point in $[1,n-1] = \{1, \ldots, n-1\}$.
Use the product rule to expand
\begin{equation}
\left( \frac{\partial}{\partial g} \right)^{n-1} \prod_{\ell \in U_n} \beta(g)  = \sum_{\phi: [n-1] \to U_n}
\prod_{\ell \in U_n     } \left( \frac{\partial}{\partial g} \right)^{\phi^{-1}(\ell)} \beta(g).
\end{equation}
The sum is over functions $\phi$ that pick out for each of the $n-1$ derivatives the factor to which it applies.

Use the Pr\"ufer correspondence. Given a linear order on $U_n$, there is a corresponding
bijection between functions $\phi$ from $[n-1]$ to $U_n$ and rooted trees with vertex
set $U_n$. Furthermore, given a rooted tree $T$, the number of times the function assumes
value $j$ in $U_n$ is the degree $ |T^{-1}(j)|$.
\end{proof}

The above may be expressed in an elegant way as
\begin{equation}
f(t,g) = \sum_{T \in \Ac_\emptyset} \frac{1}{|T|!} t^{|T|} f_T(g).
\end{equation}
Here $|T|$ denotes the number of vertices of the rooted tree, and there is one label vertex set for each
value of this number.

\begin{remark}
The problem of counting labeled rooted trees is the special case when $\beta(g) = \exp(g)$. In that case
$f_n(g) = n^{n-1} \exp(ng)$ and each $f_T(g) = \exp(ng)$. In this special case the result is
Cayley's formula $n^{n-1} = | \Ac[U_n] |$.
\end{remark}

The solution of the fixed point equation may  be expressed more economically in terms of unlabeled rooted trees by
\begin{equation}
f(t,g) = \sum_{\tau \in\tilde \Ac_\emptyset} \frac{r(\tau)}{|\tau|!} t^{|\tau|} f_\tau(g) = \sum_{\tau \in\tilde \Ac_\emptyset} \frac{1}{\sigma(\tau)} t^{|\tau|} f_\tau(g).
\end{equation}
Here $r(\tau)$ is the corresponding number of labeled rooted trees,
$|\tau|$ is the number of vertices, and $\sigma(\tau)$ is the symmetry factor.

\begin{example} Use the notation $f_\tau(g)$ to denote the factor associated with unlabeled rooted tree $\tau$. Then
\begin{eqnarray}
\lefteqn{ f(x,t) = g + f_0 t + f_1(g) t^2 + \left( \frac{1}{2} f_2(g)  + f_{1[2]}(g)  \right) t^2 } \\
&& + \left( \frac{1}{6} f_3(g) + f_{11[1]}(g)  +  \frac{1}{2} f_{1[2]}(g)  + f_{1[1[1])}(g) \right) t^4 + \cdots  . \nonumber
\end{eqnarray}
Explicitly, this is
\begin{eqnarray}
\lefteqn{ f(x,t) = g + \beta(g ) t + \beta'(g)\beta(g)  t^2 + \left( \frac{1}{2} \beta''(g)\beta(g)^2   +\beta'(g)^2\beta(g) \right) t^3 } \\
&& +
\left( \frac{1}{6} \beta'''(g)\beta(g)^3 + \beta''(g)\beta'(g) \beta(g)^2  +  \frac{1}{2}\beta''(g) \beta'(g)\beta(g)^2  + \beta'(g)^3\beta(g) \right) t^4 + \cdots .
\nonumber
\end{eqnarray}
To determine the contribution of a rooted tree, all that is needed is to know the number of vertices with given degree.
\end{example}

\section{Increasing rooted trees}

\noindent\textbf{Increasing rooted trees as partially ordered sets}
The convention used here is that the partial order of a labeled rooted tree increases as one moves away from the root.
Consider a non-empty label set $U$ together with a given linear order $\leqlin$.  An \emph{increasing  rooted tree} is a labeled rooted
tree with the property that  the map from $U$ with its
 rooted tree partial order $\leqtree$ to $U$ with the given linear order $\leqlin$ is order-preserving.
In other words, if $i \leqtree j$ in the partial order of the labeled rooted tree, then the $ i \leqlin j$ in the linear order of the labels.
For an increasing rooted tree the root $r$ is the least element in both orders. The greatest element in the linear order is maximal in
the partial order, so it is a leaf.

\noindent\textbf{Increasing rooted trees as functions}
Consider a non-empty label set $U$ together with a given linear order.  An {increasing  rooted tree} is a labeled rooted
tree $T: U \setminus \{r\} \to U$ that is  decreasing with respect to the linear order. In other words, it is
required that $T(j) \leqlin j$ in the linear order for all $j \neq r$.  The collection of all increasing rooted trees on $U$ is denoted $\Ac^\uparrow[U]$.

The relation between the three kinds of rooted trees for a given linearly ordered label set $U$ with $n\geq 1$ elements is
\begin{equation}
\Ac^\uparrow[U] \to \Ac[U] \to \tilde \Ac[n].
\end{equation}
The first map is an injection, and the second map is a surjection. The composite map is also a surjection.

\noindent\textbf{Increasing rooted trees defined recursively from below}
For each $k$ in $U$ the set of all $\ell$ with $k \leqtree \ell$ is also an increasing rooted tree $T_k$, with root $k$.
If $r$ is the root, then it must be the least element of $U$. So for each immediate successor $j$ of $r$ the rooted tree $T_j$ is an increasing
rooted tree. This gives a recursive characterization of
an increasing rooted tree on $U$  as a forest of increasing rooted trees on  $U \setminus \{r\}$.

\noindent\textbf{Increasing rooted trees defined recursively from above}
This gives another recursive description. Let $m$ be the greatest element of $U$ in the linear order.
 An increasing rooted tree on $U$ is a increasing rooted tree on $U\setminus \{m\}$ together with a point in $U \setminus \{m\}$.
 Thus there is an edge in the tree from $m$ to the chosen point. By taking $m = n, n-1, n-2, \ldots, 2$ this defines a map $\phi$ from $\{2,3, \ldots, n\}$
to the set of non-leaf vertices. (This is a variation on the Pr\"ufer  correspondence.)

\noindent\textbf{Increasing rooted trees as permutations} Each increasing rooted tree on $\{1, \ldots, n\}$ may be coded as a permutation
of $\{2, \ldots, n\}$. Such a permutation may be represented as a list of the elements of $\{2, \ldots, n\}$ in some order.
For $n=2$ the only entry in the list is 2. Suppose $n \geq 3$ and an increasing rooted trees on
$\{1, \ldots, n-1\}$  is coded as a list taken from $\{2, \ldots, n-1\}$. Consider a new increasing rooted tree on $\{1, \ldots, n\}$.
If the tree sends $n$ to $k$ with $1 \leq k \leq n-1$, then
create a new list such that for $j<k$ the $j$th place entry is the same,  the $k$th place entry is $n$, and for $j> k$ the entry in the
$j$th place is the original $j-1$ place entry.  This represents the new increasing rooted tree as a
list taken from $\{2, \ldots, n\}$.

\begin{example} The tree $1[2]$ is encoded by 2. The trees $1[23]$ and $1[2[3]]$ are encoded by $32$ and $23$. The trees $1[234]$ and $1[32[4]]$ and $ 1[23[4]]$
are encoded by $432$ and $342$ and $324$, while the trees $1[42[3]]$ and $1[2[34]] $ and $ 1[2[3[4]]]$ are encoded by $423 $ and $243$ and $ 234$.
\end{example}

The permutation representation immediately gives
the following result.

 \begin{proposition} The number of increasing rooted trees on a label set with
$n$ vertices is $(n-1)!$.
\end{proposition}

 The \emph{rooted tree factorial} $T!$ of a labeled rooted tree $T$ with root $r$ is defined inductively as the number of vertices of
$T$ times the product over $i$  with $T(i) = r$ of  $T_i!$. (An empty product gives 1.)
This is an invariant under isomorphism, so for each unlabeled rooted tree $\tau$ there is a rooted tree factorial $\tau!$.
The rooted tree factorial satisfies the recursive relation
\begin{equation}
\tau! = |\tau| \prod_{\tau'} \left( \tau'! \right)^{N(\tau')},
\end{equation}
where $N$ counts the subtrees obtained by removing the root.

There is another formula for the rooted tree factorial that is often convenient. For an unlabeled rooted tree $\tau$
consider a corresponding labeled rooted tree $T$. Let $T_j$ be the subtree over vertex $j$. Then
\begin{equation}
\tau! = \prod_j |T_j|,
\end{equation}
where the product is over all vertices of $T$. The quantity $|T_j|$ is the number of vertices of $T_j$.

The number of increasing rooted trees per unlabeled rooted tree  satisfies a recursion relation
\begin{equation}
i(\tau) = C(N) \prod_{\tau'} i(\tau')^{N(\tau')} ,
\end{equation}
where $N(\tau')$ counts immediate successor rooted trees $\tau'$.
This is similar to the formula $r(\tau)$ for the number of rooted trees per unlabeled rooted tree;
the distinction is that  there is only one choice of root point.
It follows that the ratio $r(\tau)/i(\tau)$ satisfies
\begin{equation}
\frac{r(\tau)}{i(\tau)} =  |\tau| \prod_{\tau'} \left( \frac{r(\tau')}{i(\tau')} \right)^{N(\tau')},
\end{equation}
where $N$ counts the successor rooted trees obtained by removing the root. This leads to the following relation.

\begin{proposition} Fix $n$ and an unlabeled rooted tree $\tau$ with $n$ vertices.
The ratio of the number of labeled rooted trees
to the number of increasing  rooted trees is the rooted tree factorial
\begin{equation}
\frac{r(\tau)}{i(\tau)} = \tau!.
\end{equation}
 As a consequence, the number
of increasing rooted trees for given unlabeled rooted tree $\tau$ is
\begin{equation}
i(\tau) = \frac{r(\tau)}{\tau!} = \frac{|\tau|!}{\tau! \sigma(\tau)}
\end{equation}
\end{proposition}

For $\tau$ with $n$ vertices there is an identity that expresses the fact that the sum over unlabeled rooted trees of
the corresponding number of increasing labeled rooted trees is the total number of increasing labeled rooted trees.

\begin{proposition} The sum over unlabeled rooted trees with $n$ vertices of the corresponding number of increasing rooted trees gives
\begin{equation}
\sum_{\tau \in \tilde \Ac[n]} \frac{n!}{\sigma(\tau) \tau!} = (n-1)!.
\end{equation}
\end{proposition}

\begin{example} When $n = 3$ there are only two increasing rooted trees. The symmetry factors are 2 and 1, while the
rooted tree factorials are 3 and 6. This is illustrated in Figure~8.
\end{example}

% Figure 8 Increasing rooted trees on 3 vertices
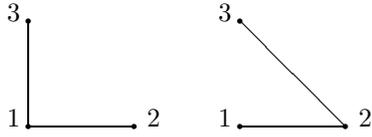
\begin{figure}[tb]
\begin{picture}(\bggl,\gl)(\negtwogl,0)
\put(0,0){
\begin{picture}(\gl,\gl)
\put(0,0){\line(1,0){\gl}}
\put(0,0){\line(0,1){\gl}}
%\put(0,0){\line(1,1){\gl}}
%\put(\gl,\gl){\line(-1,0){\gl}}
%\put(\gl,\gl){\line(0,-1){\gl}}
%\put(0,\gl){\line(1,-1){\gl}}
\put(-8,0){1}
\put(45,0){2}
\put(-8,\gl){3}
\put(0,0){\circle*{2}}
\put(0,\gl){\circle*{2}}
\put(\gl,0){\circle*{2}}
%\put(\gl,\gl){\circle*{2}}
%\put(45,\gl){4}
%\put(0,0){\circle*{4}}
%\put(\neghalfgl,\halfgl){4}
\end{picture}}
\put(\twogl,0){
\begin{picture}(\gl,\gl)
\put(0,0){\line(1,0){\gl}}
%\put(0,0){\line(0,1){\gl}}
%\put(0,0){\line(1,1){\gl}}
%\put(\gl,\gl){\line(-1,0){\gl}}
%\put(\gl,\gl){\line(0,-1){\gl}}
\put(0,\gl){\line(1,-1){\gl}}
\put(-8,0){1}
\put(45,0){2}
\put(-8,\gl){3}
%\put(45,\gl){4}
\put(0,0){\circle*{2}}
\put(0,\gl){\circle*{2}}
\put(\gl,0){\circle*{2}}
%\put(\gl,\gl){\circle*{2}}
%\put(0,0){\circle*{4}}
%\put(\neghalfgl,\halfgl){4}
\end{picture}}
\end{picture}
\caption{Increasing rooted trees on three  vertices: $\tau! = 3,6$}
\end{figure}

\begin{example}
The case $n=4$ is more interesting. There are $3!=6$ increasing rooted trees.  These map to the four unlabeled rooted trees, which
are $3, 11[1], 1[2],  1[1[1]]$.
These four unlabeled rooted trees have symmetry factors $6,1,2,1$ and rooted tree factorials $4, 8, 12, 24$.
Three of the increasing rooted trees correspond to the unlabeled rooted tree 11[1] with symmetry factor 1 and rooted tree factorial 8.
 The  number of increasing labeled rooted trees is the sum
$24/(6\cdot 4) + 24/(1 \cdot 8) + 24/(2 \cdot 12) + 24/(1\cdot 24) = 1 + 3 + 1 + 1 = 6 $. The picture is in Figure~9.
\end{example}

% Figure 9 Increasing rooted trees on 4 vertices
\begin{figure}[tb]
\begin{picture}(\bggl,\fourgl)(\neggl,0)
\put(0,\twogl){
\begin{picture}(\gl,\gl)
\put(0,0){\line(1,0){\gl}}
\put(0,0){\line(0,1){\gl}}
\put(0,0){\line(1,1){\gl}}
%\put(\gl,\gl){\line(-1,0){\gl}}
%\put(\gl,\gl){\line(0,-1){\gl}}
%\put(0,\gl){\line(1,-1){\gl}}
\put(-8,0){1}
\put(45,0){2}
\put(-8,\gl){3}
\put(45,\gl){4}
\put(0,0){\circle*{2}}
\put(0,\gl){\circle*{2}}
\put(\gl,0){\circle*{2}}
\put(\gl,\gl){\circle*{2}}
%\put(0,0){\circle*{4}}
%\put(\neghalfgl,\halfgl){4}
\end{picture}}
\put(\twogl,\twogl){
\begin{picture}(\gl,\gl)
\put(0,0){\line(1,0){\gl}}
\put(0,0){\line(0,1){\gl}}
%\put(0,0){\line(1,1){\gl}}
%\put(\gl,\gl){\line(-1,0){\gl}}
\put(\gl,\gl){\line(0,-1){\gl}}
%\put(0,\gl){\line(1,-1){\gl}}
\put(-8,0){1}
\put(45,0){2}
\put(-8,\gl){3}
\put(45,\gl){4}
\put(0,0){\circle*{2}}
\put(0,\gl){\circle*{2}}
\put(\gl,0){\circle*{2}}
\put(\gl,\gl){\circle*{2}}
%\put(0,0){\circle*{4}}
%\put(\neghalfgl,\halfgl){4}
\end{picture}}
\put(\fourgl,\twogl){
\begin{picture}(\gl,\gl)
\put(0,0){\line(1,0){\gl}}
\put(0,0){\line(0,1){\gl}}
%\put(0,0){\line(1,1){\gl}}
\put(\gl,\gl){\line(-1,0){\gl}}
%\put(\gl,\gl){\line(0,-1){\gl}}
%\put(0,\gl){\line(1,-1){\gl}}
\put(-8,0){1}
\put(45,0){2}
\put(-8,\gl){3}
\put(45,\gl){4}
\put(0,0){\circle*{2}}
\put(0,\gl){\circle*{2}}
\put(\gl,0){\circle*{2}}
\put(\gl,\gl){\circle*{2}}
%\put(0,0){\circle*{4}}
%\put(\neghalfgl,\halfgl){4}
\end{picture}}
\put(0,0){
\begin{picture}(\gl,\gl)
\put(0,0){\line(1,0){\gl}}
%\put(0,0){\line(0,1){\gl}}
\put(0,0){\line(1,1){\gl}}
%\put(\gl,\gl){\line(-1,0){\gl}}
%\put(\gl,\gl){\line(0,-1){\gl}}
\put(0,\gl){\line(1,-1){\gl}}
\put(-8,0){1}
\put(45,0){2}
\put(-8,\gl){3}
\put(45,\gl){4}
\put(0,0){\circle*{2}}
\put(0,\gl){\circle*{2}}
\put(\gl,0){\circle*{2}}
\put(\gl,\gl){\circle*{2}}
%\put(0,0){\circle*{4}}
%\put(\neghalfgl,\halfgl){4}
\end{picture}}
\put(\twogl,0){
\begin{picture}(\gl,\gl)
\put(0,0){\line(1,0){\gl}}
%\put(0,0){\line(0,1){\gl}}
%\put(0,0){\line(1,1){\gl}}
%\put(\gl,\gl){\line(-1,0){\gl}}
\put(\gl,\gl){\line(0,-1){\gl}}
\put(0,\gl){\line(1,-1){\gl}}
\put(-8,0){1}
\put(45,0){2}
\put(-8,\gl){3}
\put(45,\gl){4}
\put(0,0){\circle*{2}}
\put(0,\gl){\circle*{2}}
\put(\gl,0){\circle*{2}}
\put(\gl,\gl){\circle*{2}}
%\put(0,0){\circle*{4}}
%\put(\neghalfgl,\halfgl){4}
\end{picture}}
\put(\fourgl,0){
\begin{picture}(\gl,\gl)
\put(0,0){\line(1,0){\gl}}
%\put(0,0){\line(0,1){\gl}}
%\put(0,0){\line(1,1){\gl}}
\put(\gl,\gl){\line(-1,0){\gl}}
%\put(\gl,\gl){\line(0,-1){\gl}}
\put(0,\gl){\line(1,-1){\gl}}
\put(-8,0){1}
\put(45,0){2}
\put(-8,\gl){3}
\put(45,\gl){4}
\put(0,0){\circle*{2}}
\put(0,\gl){\circle*{2}}
\put(\gl,0){\circle*{2}}
\put(\gl,\gl){\circle*{2}}
%\put(0,0){\circle*{4}}
%\put(\neghalfgl,\halfgl){4}
\end{picture}}
\end{picture}
\caption{Increasing rooted trees on four  vertices: $\tau! = 4, 8, 8, 12, 24$}
\end{figure}
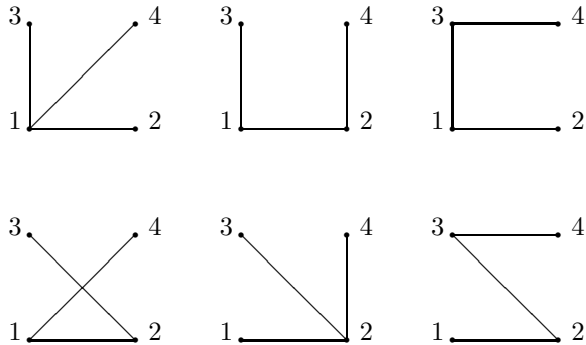

\begin{example}
For $n=5$ vertices there are 9 unlabeled rooted trees. These are indicated in Table~1. The symmetry group of all permutations has
order $n!= 120$. The table lists the  symmetry
factors $\sigma(\tau)$, the number of labeled rooted trees $r(\tau)$, the tree factorial $\tau!$, and the number
of unlabeled rooted trees. The symmetry factor $\sigma(\tau)$ may be read off from the multiset description of $\tau$. The other quantities
are related by $r(\tau) = n!/\sigma(\tau)$ and $i(\tau) = r(\tau)/\tau!$. If all is well, the total number of labeled rooted trees should be $n^{n-1} = 625$, while the total number of increasing rooted trees should be $(n-1)! = 24$.
\end{example}

\begin{remark}
It is instructive to look at the exponential generating function for the number $a^\uparrow_n$ of increasing rooted trees
with $n$ vertices.
The recursive definition  gives
\begin{equation}
\frac{d}{dt} a^\uparrow(t) =  \exp( a^\uparrow(t)) .
\end{equation}
This has the easy solution $a^\dagger(t) = \log( \frac{1}{1-t} )$.
The details are in \cite{BergeronLabelleLeroux}.
\end{remark}

\begin{table}
\begin{center}
\begin{tabular}{lrrrr}
$\tau$ &   $\sigma(\tau)$  & $r(\tau)$ & $\tau!$  & $i(\tau)$   \\ \hline
4 &  24 & 5 & 5 & 1 \\
21[1] & 2 & 60 & 10 & 6 \\
2[1] & 2 & 60 & 20 & 3 \\
11[2] &2 & 60 & 15 & 4 \\
11[1[1]] &1 & 120 & 30 & 4 \\
1[3]   &6 & 20 & 20 & 1 \\
1[11[1]] &1 & 120 & 40 &3 \\
1[1[2]] &2& 60 & 60 & 1 \\
1[1[1[1]]] &1 & 120 & 120 & 1 \\
\end{tabular}
\caption{Unlabeled rooted trees on 5 vertices}
\end{center}
\end{table}

\section{Ordinary differential equations and increasing rooted trees}

Let $\beta(x)$ be a formal power series in $x$. Let $t$ and $g$ be parameters. Consider
the \emph{ordinary differential equation}
\begin{equation}
\frac{dx}{dt} = \beta(x)
\end{equation}
with initial condition $x = g$ at $t = 0$.

\begin{proposition} This ordinary differential equation  has the formal solution
\begin{equation}
x = \bar f(t,g) = \sum_{n=0}^\infty \frac{t^n}{n!} \bar f_n(g) ,
\end{equation}
where  the coefficient $\bar f_n(g)$ has the explicit representation
\begin{equation}
\bar f_n(g) = \left( \beta(g) \frac{\partial}{\partial g} \right)^n g.
\end{equation}
\end{proposition}

\begin{proof} Let the solution of the initial value problem be $\bar f(t,g)$.
It is easy to see by induction that
\begin{equation}
\frac{\partial^n \bar f(t,g)}{\partial t^n} = h_n(\bar f(t,g))
\end{equation}
for a suitable function $h_n(x)$.
Taking one more derivative gives
\begin{equation}
h_{n+1}(\bar f(t,g))  = h'_n(\bar f(t,g)) \beta(\bar f(t,g)) .
\end{equation}
Setting $t=0$ gives $\bar f(0,g) = g$, so the result is
\begin{equation}
h_{n+1}(g) = h'_n(g) \beta(g) .
\end{equation}
The conclusion follows immediately.
\end{proof}

The combinatorial solution is given by an expansion indexed by increasing rooted trees.
For each $n\geq 1$  fix a linearly ordered label set $U_n$ with $k$ elements. For instance,  take $U_n = \{1, \ldots, n \}$ with
the usual linear order.
Furthermore,  consider increasing rooted trees $\bar T$ with the label set as vertex set.
For $n = 0$ introduce an empty set object.

\begin{proposition} The ordinary differential equation  has the  solution
given by a formal power series as above,
where  $\bar f_n(g)$ has the explicit representation
\begin{equation}
\bar f_n(g) = \sum_{\bar T \in \Ac^\uparrow [ U_n]} f_{\bar T}(g)
\end{equation}
For $n \geq 1$ the coefficient is
\begin{equation}
 f_{\bar T}(g) = \prod_{j \in U_n} \left( \frac{\partial}{\partial g} \right)^{|\bar T^{-1}(j)| } \beta(g),
\end{equation}
where $|{\bar T}^{-1}(j)|$ is the degree of vertex $j$ of rooted tree ${\bar T}$.
For $n=0$ the contribution from
the empty rooted tree is $f_\emptyset(g) = g$.
\end{proposition}

\begin{proof}
Write
\begin{equation}
\bar f_n(g) = \left( \beta(g) \frac{\partial}{\partial g} \right)^{n-1} \beta(g).
\end{equation}
Index the partial derivatives from $n$ down to $2$. Index the $\beta(g)$ factors from $n$ down to 1.
Then every partial derivative acts only on the $\beta(g)$ factors with strictly smaller index.
So
\begin{equation}
\bar f_n(g) = \sum_\phi \prod_{j=1}^n \left( \frac{\partial}{\partial g} \right)^{|\phi^{-1}(j)|} \beta(g),
\end{equation}
where the sum is over functions $\phi$ from to $[2,n]$ to $[1,n]$ with the property
that $\phi(i) < i$ for all $i$.

Every such function $\phi$ is an increasing rooted tree  on
$[1,n]$, where $\phi(i)$ is the immediate predecessor of $i$, and $\phi^{-1}(j)$ is
the set of immediate successors of $j$.
\end{proof}

The above may be expressed in an elegant way as
\begin{equation}
\bar f(t,g) = \sum_{\bar T \in \Ac_\emptyset^\uparrow}  \frac{1}{|\bar T|!} t^{|\bar T|}  f_{\bar T}(g).
\end{equation}
Here $|\bar T|$ denotes the number of vertices of the increasing rooted tree, and there is one label vertex set for each
value of this number.

\begin{remark}
The problem of counting increasing rooted trees is the special case when $\beta(g) = \exp(g)$. The solution of the
differential equation is given explicitly by $x = g - \log(1-\exp(g)t)$. The coefficient in $n$th order is
$\bar f_n(g) = (n-1)! \exp(ng)$, and each $f_T(g) = \exp(ng)$. In this special case the result is equivalent to
the formula  $(n-1)! = | \Ac^\uparrow [U_n] |$.
\end{remark}

The solution of the ordinary differential equation
may be expressed in terms of unlabeled rooted trees by
\begin{equation}
\bar f(t,g) = \sum_{\tau \in \tilde \Ac_\emptyset} \frac{i(\tau)}{|\tau|!} t^{|\tau|} f_\tau(g) = \sum_{\tau \in \tilde \Ac_\emptyset} \frac{1}{\sigma(\tau) \tau!} t^{|\tau|} f_\tau(g).
\end{equation}
Here $r(\tau)$ is the corresponding number of increasing rooted trees,
$|\tau|$ denotes the number of vertices of the rooted tree, $\sigma(\tau)$ is the symmetry factor, and $\tau!$ is the
rooted tree factorial.

\begin{example} Use the notation $f_\tau(g)$ to denote the factor associated with unlabeled rooted tree $\tau$. Then
\begin{eqnarray}
\lefteqn{ \bar f(x,t) = g + f_0 t + \frac{1}{2} f_1(g) t^2 + \frac{1}{6} \left(  f_2(g)  + f_{1[2]}(g)  \right) t^2 } \\
&& + \frac{1}{24} \left(  f_3(g) + 3 f_{11[1]}(g)  +   f_{1[2]}(g)  + f_{1[1[1]]}(g) \right) t^4 + \cdots  . \nonumber
\end{eqnarray}
Explicitly, this is
\begin{eqnarray}
\lefteqn{ \bar f(x,t) = g + \beta(g) t + \frac{1}{2}\beta'(g)\beta(g)  t^2 + \frac{1}{6} \left(  \beta''(g)\beta(g)^2   + \beta'(g)^2\beta(g) \right) t^3 } \\
&& +
\frac{1}{24} \left(  \beta'''(g)\beta(g)^3 + 3\beta''(g)\beta'(g) \beta(g)^2  +  \beta''(g) \beta'(g)\beta(g)^2  + \beta'(g)^3\beta(g) \right) t^4 + \cdots .
\nonumber
\end{eqnarray}
The factor 3 comes from the 3 increasing rooted trees associated with $\tau = 11(1)$. These ordinary differential equation coefficients
are the fixed point equation coefficients divided by tree factorials.
\end{example}

\section{The Butcher group (composition) for labeled rooted trees}

Let $T$ be a labeled rooted tree with vertex set $U$ and root $r$ in $U$.
 A \emph{rooted subtree} $T_0$ is a rooted tree on  some non-empty subset $U_0 \subseteq U$ with the same root $r$ in $U_0$ that is a restriction of $T$ to this subset. The condition that $T_0$ is a rooted subtree of $T$ is denoted $T_0 \to T$. The
 empty subset $U_0 = \emptyset$ corresponds to an empty set object $T_0$; that case is also abbreviated $T_0 \to T$.

 For each  $T_0$ with $T_0 \to T$ there is a corresponding \emph{difference forest} $T \setminus T_0$ of rooted trees on $U \setminus U_0$.
 The trees in the forest are all subtrees $T_j$ with $j \in U \setminus U_0$ and $T(j) \in U_0$. If $U_0 = \emptyset$ and $T_0$ is the empty set object, then
 the difference forest consists of the tree $T$ on $U$.

\begin{proposition}
For each non-empty label set $U$ there is a one-to-one correspondence between rooted tree pairs $T_0, T$ with $T_0 \to T$ and triples $T_0$, $F_1$, $\phi$, where $T_0$ is a rooted tree on $U_0 \subseteq U$, $F_1$ is a forest on $U_1 = U \setminus U_0$, and
$\phi$ is a function from the set partition $\Gamma_1$ of the forest to $U_0$.
\end{proposition}

Let $\tilde \Ac_\emptyset$ be the set of unlabeled rooted trees together with the empty set object.
Consider the space $\Rb^{\tilde \Ac_\emptyset}$ of all functions $c$ from $\tilde \Ac_\emptyset$ to the real numbers.
These are coefficients $c(\tau)$
that depend on unlabeled rooted trees $\tau$. Since each labeled rooted tree $T$ determines a corresponding unlabeled rooted tree $\tau$,
the coefficients $c(T)$ are also defined for labeled rooted trees.
For a forest of rooted trees the coefficient $a^\times(F)$ is the product of the $a(T')$ for $T'$ in the forest. In particular, for the
empty forest $a^\times(\emptyset) =1$.

The operation of \emph{subtree convolution} $a * b$ is defined  when $a(\emptyset) = 1$  by $c = a*b$, where
\begin{equation}
c(T) = \sum_{T_0: T_0 \rightarrow T} b(T_0) a^\times(T \setminus T_0).
\end{equation}
When $T_0$ is the empty set object, the corresponding term in the sum is is $b(\emptyset) a(T)$. When  $T_0 = T$,
the corresponding term in the sum is $b(T)a^\times(\emptyset) = b(T)$.
When $T$ is the empty set object, then $c(\emptyset) = b(\emptyset)$.

In the multiplication $a*b$ the forest factor is on the left and the tree factor is on the right.
This convention is common in this context \cite{HairerWanner2,ChartierHairerVilmart}.
The Butcher group multiplication is the special case when both $a(\emptyset) = 1$ and $b(\emptyset) = 1$. This group is denoted $G_C$, where $C$ stands for
composition.  The identity $\delta_\emptyset$ in the group has coefficient 1 for the empty set object and 0 for the rooted trees.
An inductive argument shows that every element has an inverse element. The group $G_C$ is the character group of the rooted tree Hopf algebra
of Connes and Kreimer.

\begin{example} It is easy to compute the Butcher group multiplication for reasonably small labeled rooted trees. Here are the results for up to 3 vertices.
For notational simplicity label the tree as an increasing tree. Take the root with label 1. For one-vertex rooted trees  $c(1) = b(1) + a(1)$. For two-vertex rooted trees the result is $c(1[2]) = b(1[2]) + b(1)a(2) + a(1[2])$. The first
interesting case is $n=3$. Let $1[23]$ be the rooted tree with root at 1 and with leaves at 2 and 3, and let $1[2[3]]$ be the rooted tree with root at 1 and with successor rooted tree $2[3]$. Then
\begin{equation}
c(1[23]) = b(1[23]) + b(1[2])a(3) + b(1[3])a(2) + b(1)a(2)a(3) + a(1[23]).
\end{equation}
Similarly
\begin{equation}
c(1[2[3]]) = b(1[2[3]]) + b(1[2])a(3)  + b(1)a(2[3]) + a(1[2[3]]).
\end{equation}
The expressions for the two cases are quite different.
\end{example}

\begin{example}
Here is one calculation for $n=4$. The rooted tree is $T = 1[3[2[4]]$ with root at 1 and successor rooted trees 3 and $2[4]$. This happens to be
an increasing rooted tree, but that is just for notational convenience. The result is
\begin{eqnarray}
\lefteqn{c(T) = b(T) + b(1[2[4]])a(3) + b(1[23])a(4)} \kill \\
&&  + b(1[3])a(2[4]) + b(1[2])a(3)a(4)) + b(1)a(2[4])a(3)  + a(T). \nonumber
\end{eqnarray}
This is illustrated in Figure~10. The number of vertices in the subtrees are 4, 3, 3, 2, 2, 1, 0. The difference
forests consist of 0,  1, 1, 1, 2, 2, 1 rooted trees.
\end{example}

% Figure 10 Subtrees of a rooted tree:
\begin{figure}[tb]
\begin{picture}(\bggl,\fourgl)(\neggl,0)
\put(0,\twogl){
\begin{picture}(\gl,\gl)
\put(0,0){\line(1,0){\gl}}
\put(0,0){\line(0,1){\gl}}
%\put(0,0){\line(1,1){\gl}}
%\put(\gl,\gl){\line(-1,0){\gl}}
\put(\gl,\gl){\line(0,-1){\gl}}
%\put(0,\gl){\line(1,-1){\gl}}
%\put(-8,0){1}
%\put(45,0){2}
%\put(-8,\gl){3}
%\put(45,\gl){4}
\put(0,0){\circle*{6}}
\put(0,\gl){\circle*{2}}
\put(\gl,0){\circle*{2}}
\put(\gl,\gl){\circle*{2}}
%\put(0,0){\circle*{4}}
%\put(\neghalfgl,\halfgl){4}
\end{picture}}
\put(\twogl,\twogl){
\begin{picture}(\gl,\gl)
\put(0,0){\line(1,0){\gl}}
%\put(0,0){\line(0,1){\gl}}
%\put(0,0){\line(1,1){\gl}}
%\put(\gl,\gl){\line(-1,0){\gl}}
\put(\gl,\gl){\line(0,-1){\gl}}
%\put(0,\gl){\line(1,-1){\gl}}
%\put(-8,0){1}
%\put(45,0){2}
%%\put(-8,\gl){3}
%\put(45,\gl){4}
\put(0,0){\circle*{6}}
\put(0,\gl){\circle*{4}}
\put(\gl,0){\circle*{2}}
\put(\gl,\gl){\circle*{2}}
%\put(0,0){\circle*{4}}
%\put(\neghalfgl,\halfgl){4}
\end{picture}}
\put(\fourgl,\twogl){
\begin{picture}(\gl,\gl)
\put(0,0){\line(1,0){\gl}}
\put(0,0){\line(0,1){\gl}}
%\put(0,0){\line(1,1){\gl}}
%\put(\gl,\gl){\line(-1,0){\gl}}
%\put(\gl,\gl){\line(0,-1){\gl}}
%\put(0,\gl){\line(1,-1){\gl}}
%\put(-8,0){1}
%\put(45,0){2}
%\put(-8,\gl){3}
%\put(45,\gl){4}
\put(0,0){\circle*{6}}
\put(0,\gl){\circle*{2}}
\put(\gl,0){\circle*{2}}
\put(\gl,\gl){\circle*{4}}
%\put(0,0){\circle*{4}}
%\put(\neghalfgl,\halfgl){4}
\end{picture}}
\put(0,0){
\begin{picture}(\gl,\gl)
%\put(0,0){\line(1,0){\gl}}
\put(0,0){\line(0,1){\gl}}
%\put(0,0){\line(1,1){\gl}}
%\put(\gl,\gl){\line(-1,0){\gl}}
\put(\gl,\gl){\line(0,-1){\gl}}
%\put(0,\gl){\line(1,-1){\gl}}
%\put(-8,0){1}
%\put(45,0){2}
%\put(-8,\gl){3}
%\put(45,\gl){4}
\put(0,0){\circle*{6}}
\put(0,\gl){\circle*{2}}
\put(\gl,0){\circle*{4}}
\put(\gl,\gl){\circle*{2}}
%\put(0,0){\circle*{4}}
%\put(\neghalfgl,\halfgl){4}
\end{picture}}
\put(\twogl,0){
\begin{picture}(\gl,\gl)
\put(0,0){\line(1,0){\gl}}
%\put(0,0){\line(0,1){\gl}}
%\put(0,0){\line(1,1){\gl}}
%\put(\gl,\gl){\line(-1,0){\gl}}
%\put(\gl,\gl){\line(0,-1){\gl}}
%\put(0,\gl){\line(1,-1){\gl}}
%\put(-8,0){1}
%\put(45,0){2}
%\put(-8,\gl){3}
%\put(45,\gl){4}
\put(0,0){\circle*{6}}
\put(0,\gl){\circle*{4}}
\put(\gl,0){\circle*{2}}
\put(\gl,\gl){\circle*{4}}
%\put(0,0){\circle*{4}}
%\put(\neghalfgl,\halfgl){4}
\end{picture}}
\put(\fourgl,0){
\begin{picture}(\gl,\gl)
%\put(0,0){\line(1,0){\gl}}
%\put(0,0){\line(0,1){\gl}}
%\put(0,0){\line(1,1){\gl}}
%\put(\gl,\gl){\line(-1,0){\gl}}
\put(\gl,\gl){\line(0,-1){\gl}}
%\put(0,\gl){\line(1,-1){\gl}}
%\put(-8,0){1}
%\put(45,0){2}
%\put(-8,\gl){3}
%\put(45,\gl){4}
\put(0,0){\circle*{6}}
\put(0,\gl){\circle*{4}}
\put(\gl,0){\circle*{4}}
\put(\gl,\gl){\circle*{2}}
%\put(0,0){\circle*{4}}
%\put(\neghalfgl,\halfgl){4}
\end{picture}}
\put(\sixgl,0){
\begin{picture}(\gl,\gl)
\put(0,0){\line(1,0){\gl}}
\put(0,0){\line(0,1){\gl}}
%\put(0,0){\line(1,1){\gl}}
%\put(\gl,\gl){\line(-1,0){\gl}}
\put(\gl,\gl){\line(0,-1){\gl}}
%\put(0,\gl){\line(1,-1){\gl}}
%\put(-8,0){1}
%\put(45,0){2}
%\put(-8,\gl){3}
%\put(45,\gl){4}
\put(0,0){\circle*{4}}
\put(0,\gl){\circle*{2}}
\put(\gl,0){\circle*{2}}
\put(\gl,\gl){\circle*{2}}
%\put(0,0){\circle*{4}}
%\put(\neghalfgl,\halfgl){4}
\end{picture}}
\end{picture}
\caption{Subtrees and difference forests}
\end{figure}
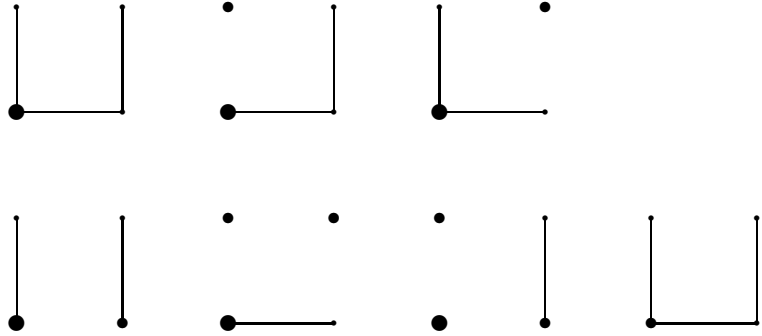

Coefficients depending on rooted trees also define certain functions (defined as formal power series).
Each such function is a weighted tree sum in the form of an exponential generating
function. For coefficients $c$ the function is
\begin{equation}
f^c(t,g) = \sum_{T \in \Ac_\emptyset} \frac{1}{|T|!} t^{|T|} c(T) f_T(g).
\end{equation}
Again for rooted tree $T$ the coefficient is
\begin{equation}
f_T(g) = \prod_{\ell \in U}\left( \frac{\partial}{\partial g} \right)^{ |T^{-1}(\ell)|} \beta(g) .
\end{equation}
The zero order term is $c(\emptyset) g$.   The remainder of the series is a power
series in powers of $t$ and powers of derivatives of $\beta(g)$. Such as series may be further
expanded in powers of $g$, but that is not done here.

\begin{theorem}[Composition] Suppose $a(\emptyset) = 1$, and let
 $c = a*b$ be the subtree convolution.
The corresponding weighted sums are related by
\begin{equation}
f^c(t,g) = f^b(t,f^a(t,g)).
\end{equation}
\end{theorem}

\begin{proof} The proof is based on the calculus formulas in the appendix.
The task is to show that $f^c(t,g)$ and $f^b(t,f^a(t,g))$ have the same terms of each order.
The $n$th term is the $n$th partial derivative with respect to $t$, evaluated at zero. (The case $n=0$ is
trivial, so  assume $n \geq 1$.) Take a set $U$ with $|U| = n$ elements. For $f^c(t,g)$
\begin{equation}
D_t^{|U|} f^c(t,g)|_{t=0}  =     \sum_{T \in \Ac[U]} c(T) f_T(g).
\end{equation}

\noindent\textbf{Product rule} The composition
\begin{equation}
f^b(t,f^a(t,g)) =   \sum_{T\in \Ac_\emptyset} \frac{1}{|T|!} t^{|T|} b(T) f_{T}( f^a(t,g)).
\end{equation}
is a combination of products of two factors $t^{|T|}$ and $f_{T}(f^a(t,g))$
By the product rule, the $n$th derivative is
a sum over disjoint unions $U = U_0 + U_1$ of the form
\begin{equation}
D_t^{|U|} f^b(t,f^a(t,g))|_{t=0}  =    \sum_{U = U_0 + U_1} \sum_{T_0 \in \Ac[U_0]} b(T_0)
 D_t^{|U_1|} f_{T_0}( f^a(t,g))|_{t=0} .
\end{equation}

\noindent\textbf{Chain rule}
The next task is to evaluate  $D^{|U_1|} f_{T_0}( f^a(t,g))$ by the chain rule.  This is a sum over set partitions $\Gamma_1$ of
$U_1$
\begin{equation}
 D_t^{|U_1|} f_{T_0}( f^a(t,g))|_{t=0} =
 \sum_{\Gamma_1 \in \setpartition[U_1]} D_g^{|\Gamma_1|} f_{T_0}(f^a(t,g))|_{t = 0}  \prod_{B \in \Gamma_1} D^{|B|} f^a(t,g)|_{t=0}.
 \end{equation}
 Since $f^a(0,g) = g$, this is
 \begin{equation}
 D_t^{|U_1|} f_{T_0}( f^a(t,g))|_{t=0} =
 \sum_{\Gamma_1 \in \setpartition[U_1]} D_g^{|\Gamma_1|} f_{T_0}(g) \prod_{B\in \Gamma_1} \left( \sum_{H \in \Ac[ B] } a(H) f_H(g) \right).
 \end{equation}

\noindent\textbf{Distributive law}
The distributive law is used to expand the product
\begin{equation}
\prod_{B\in \Gamma_1} \left( \sum_{H \in \Ac[B]} a(H) f_H(g) \right) = \sum_F \prod_{B \in \Gamma_1} a(F(B)) f_{F(B)}(g),
\end{equation}
where $F$ is a function defined on $\Gamma_1$ such that the value of $F$ on block $B$ in $\Gamma_1$ is a rooted tree on $B$.

\noindent\textbf{Product rule}
When $T_0$ is non-empty the coefficient $f_{T_0}(g)$ is a product over vertices $\ell$ in $U_0$ of factors $\beta^{(d(\ell))}(g)$.
By the product rule the derivative of order $|\Gamma_1|$ is  a sum over functions $\phi: \Gamma_1 \to U_0$ in the form
\begin{equation}
D_g^{|\Gamma_1|} f_{T_0}(g) = \sum_{\phi: \Gamma_1 \to U_0} f_{T_0,\phi}(g).
\end{equation}
Here
\begin{equation}
f_{T_0,\phi}(g) = \prod_{\ell \in U_0} \beta^{(|T_0^{-1}(\ell)| + |\phi^{-1}(\ell)|)} (g).
\end{equation}
For the empty rooted tree  $f_\emptyset(g) = g$ and so the only derivative that is non-zero is the first derivative, corresponding
to a set partition of $U = U_1$ into one block. It is convenient to define $f_{\emptyset,\phi}(g) = 1$, where $\phi$ is some unique
object whose nature is not important.

\noindent\textbf{Construction of rooted tree and subtree}
The end result is
\begin{eqnarray}
\lefteqn{D_t^{|U|} f^b(t,f^a(t,g))|_{t=0} = } \\
 && \sum_{U = U_0 + U_1} \sum_{T_0 \in \Ac[ U_0]}
\sum_{\Gamma_1 \in \setpartition[U_1]} \sum_{\phi: \Gamma_1 \to U_0}  \sum_F \prod_{B \in \Gamma_1}  b(T_0)a(F(B)) f_{T_0,\phi}(g) f_{F(B)}(g) . \nonumber
\end{eqnarray}
Fix $U = U_0 + U_1$ and rooted tree $T_0$ on $U_0$. The remaining data $\Gamma_1$, $\phi$, $F$ determine a rooted tree
$T$ on $U$ that extends $T_0$. The value $\phi(B)$ may be thought of as the vertex to which the root of the tree $F(B)$ maps.
So
\begin{equation}
D_t^{|U|} f^b(t,f^a(t,g))|_{t=0} =    \sum_{U = U_0 + U_1} \sum_{T_0 \in \Ac[ U_0] } \sum_{T: T_0 \rightarrow T}
 b(T_0)a^\times(T \setminus T_0) f_T(g) .
\end{equation}
The sum may be done in the other order, first the rooted tree $T$ and then the subtree $T_0$. This gives
\begin{equation}
D_t^{|U|} f^b(t,f^a(t,g)) |_{t=0} =     \sum_{T \in \Ac[U]} \sum_{T_0 : T_0 \rightarrow T}
 b(T_0)a^\times(T \setminus T_0) f_T(g) .
\end{equation}
In other words, the derivative is $\sum_{T \in \Ac[U]} (a * b)(T) f_T(g)$.
\end{proof}

One special case of the multiplication law is when the sequence $b$ is zero except for $b(\bullet) = 1$. This corresponds to
the function $f^b(t,g) = t \beta(g)$. In this case $c(T) = a^\times( T \setminus \bullet)$, the product of $a(T')$ for all rooted trees
$T'$ in the successor forest.

\begin{example} Consider $f^b(t,g) = t \beta(g)$ and $f^a(t,g) = g + a(\bullet) t \beta(g)$. The composition is
\begin{equation}
f^c(t,g) = t \beta( g + a(\bullet) t \beta(g)).
\end{equation}
This type of composition occurs in numerical methods for the solution
of ordinary differential equations. In this case $c(T)$ is zero unless $T$ is a rooted tree on a set with $n\geq 1$ points,
one root and $n-1$ leaves. There are $n$ such rooted trees. The corresponding successor forests each consist of $n-1$ one point rooted trees.
This gives $c(T) = a(\bullet)^{n-1}$. The function $f^c(t,g)$ has the rooted tree expansion
\begin{equation}
f^c(t,g) = \sum_{n=1}^\infty n \frac{t^n}{n!} a(\bullet)^{n-1} \beta^{(n-1)}(g) \beta(g)^{n-1}.
\end{equation}
This is the Taylor expansion of the composite function.

This example underpins the  Runge-Kutta methods for the numerical solution of ordinary differential equations.
The first order Euler method is to use $g + \beta(g) t$ to approximate the solution. Various second order methods depend on a parameter $a$;
they are of the form $ g + (1 - \frac{1}{2a} ) \beta(g) t  + \frac{1}{2a}    \beta(g + a \beta(g) t )$. In particular, $a = 1$ is
the analog of the trapezoidal rule, and $a = \frac{1}{2} $ is the analog of the midpoint rule. Since $\beta(g + a\beta(g)t)$
agrees with $\beta(g) +  a \beta'(g)  \beta(t) t $ to second order, the second order Runge-Kutta method agrees with the Taylor method
$g + \beta(t) t + \frac{1}{2} \beta'(g) \beta(g) t $ to second order. The Runge-Kutta method
 has the advantage that it does not require computing the derivative $\beta'(g)$.
\end{example}

The Butcher group is related to composition of power series. Take the case when $\beta(g) = \exp(g)$. In that case the contribution of
a labeled rooted tree on $n \geq 1$ vertices is $\exp(ng)$. Suppose $c = a * b$ and $f^c(t, g) = f^b( t, f^a(t,g)) $. Each of the individual functions is a power series in powers of $\exp(g)$.  Define $h^a(w) = \exp( f^a(1,\log(w)) )$ and similarly for the others. The resulting functions are
power series in $w$, and they are related by $h^c(w) = h^b(h^a(w))$. Subtree convolution is mapped to composition of power series.

\section{The Butcher group for increasing  rooted trees}

The composition law for the Butcher group may also be written in terms of increasing rooted trees.
The function corresponding to sequence $a$ is
\begin{equation}
f^a(t,g) = \sum_{\bar T \in \Ac_\emptyset^\uparrow} \frac{t^{|{\bar T}|}}{|{\bar T}|!} a({\bar T}) \bar T! f_{\bar T}(g).
\end{equation}
The rooted subtree convolution is defined in the same way, since if $T$ is an increasing rooted tree and $T_0 \to T$,
then the subtree $T_0$ is also increasing.

If $T_0 \to T$, the \emph{rooted tree binomial coefficient} is
\begin{equation}
 {T \choose T_0} = \frac{T!}{T_0! \prod_{T' \in T \setminus T_0} T'!} .
 \end{equation}
 For a linear rooted tree this is the usual binomial coefficient.

The change of variable $\bar c_{ T} = c({ T}) T!$ gives another representation of the Butcher multiplication
as
\begin{equation}
\bar c(T) = \sum_{T_0: T_0 \rightarrow T} { T \choose  T_0} \bar b(T_0) \bar a^\times(T \setminus T_0).
\end{equation}

The solution $\bar f(t,g)$ of the differential equation $dx/dt = \beta(x)$ with initial condition $g$ is the
case corresponding to $a( T) = 1/ T!$, or $\bar a( T) = 1$.
 This fact leads to an identity for a  binomial
coefficient associated with rooted trees.

\begin{proposition}
 For a labeled rooted tree with $n$ vertices
\begin{equation}
 \sum_{T_0 \to T}  {T \choose T_0} = 2^n.
\end{equation}
\end{proposition}

\begin{proof} The solution of the differential equation has the group property
 $\bar f(2t,g) = \bar f(t,\bar f(t,g)) $. Take $\bar a(T) = 1$ and the corresponding $\bar f^a(t,g)$.
 Then $\bar f^a(2t,g) = \bar f^a(t,\bar f^a(t,g))$. Now take
 $\bar a(T) = \bar b(T) = 1$, so $\bar c(T)=  2^{|T|}$. This then translates to $\bar f^c(T) = \bar f^b(t, \bar f^a(t,g))$.
 This gives the rooted tree binomial coefficient identity for increasing rooted trees. Since rooted tree factorials do not
 depend on the order on the label set, the identity holds for all labeled rooted trees.
\end{proof}

\begin{example}
For the labeled rooted tree with four vertices shown in Figure~10 the  binomial identity says
\[
16 = 1+ \frac{4}{3}  + \frac{8}{3}  + 2  +
4  +  4   +  1 .
\]
The rooted tree binomial coefficients need not be whole numbers, but the sum is always a  power of 2.
\end{example}

\section{The Butcher group for unlabeled rooted trees}

The Butcher group may also be presented using unlabeled rooted trees, but there is a complication.
If $\tau, \tau_0$ are unlabeled rooted trees, choose  $T$ to be a labeled rooted tree that determines $\tau$.
 Consider the set of labeled rooted trees $T_0$ with $T_0 \to T$ and with $T_0$ determining $\tau_0$.
 This set depends on the chosen $T$, but the number of elements in the set is independent of $T$. Denote this
 number by $[\tau,\tau_0]$.

The multiplication operation for coefficients may then be written $c = a * b$, where
\begin{equation}
c(\tau) = \sum_{\tau_0} [\tau, \tau_0]  b(\tau_0) a^\times(\tau \setminus \tau_0).
\end{equation}
The extra complication is the presence of the multiplicity coefficients $[\tau, \tau_0]$.

For an arbitrary coefficient function $c$ there is a corresponding function
\begin{equation}
f^c(t,g) = \sum_{\tau \in \tilde \Ac_\emptyset} \frac{1}{\sigma(\tau)} c(\tau) t^{|\tau|} f_\tau(g).
\end{equation}
The zero order term is $c(\emptyset) g$, and for the other terms
\begin{equation}
f_\tau(g) = \prod_k (D^k_g \beta(g))^{v_k(\tau)}.
\end{equation}
 This is precisely the same sum as before. The following is a restatement of the previous theorem.

\begin{corollary}
Suppose $a(\emptyset) = 1$.
Define $c = a * b$ with the multiplicity factor.
The corresponding functions
satisfy
\begin{equation}
f^c(t,g) = f^b(t,f^a(t,g)).
\end{equation}
\end{corollary}

\begin{example}
Take $k=4$ vertices, and let $3$ be the unlabeled rooted tree with a root and three leaves. Let $2$ be the unlabeled  rooted
tree with a root and two leaves.
Then the multiplicity $[3, 2] = 3$. Similarly, $[3,1] = 3$. On the other hand, $[3,0] = 1$, since there is only one way of inserting the root.
The conclusion is that
\begin{equation}
c(3) = b(3) + 3 b(2) a(0) + 3 b(1) a(0)^2 + b(0) a(0)^3 + a(3).
\end{equation}
In the world of unlabeled rooted trees, multiplicity factors are inescapable.
\end{example}

\section{Substitution for labeled rooted trees}

The Butcher group is about subtree convolution and composition of functions. There is another algebraic structure
based on quotient tree convolution and substitution. The reference \cite{ChartierHairerVilmart} has an account of
the subject and its history, along with some applications. See also \cite{CalaqueEbrahimiFardManchon} for a Hopf algebra approach and for more background.
The following is a brief account in the labeled rooted tree framework.

Given a rooted tree $T$ on $U$ with root $r$ and a set $R$ with $r \in R$, there is a corresponding forest $F$ obtained by restricting
$T$ to $U \setminus R$. This is called a \emph{subforest} of the rooted tree. The set partition $\Gamma$ defined by this forest is in
one-to-one correspondence with the set of roots $R$ of the rooted trees in the forest. When $F$ is a subforest of $T$ we can write
$F \sqsubseteq T$.

Given a rooted tree $T$ and a subforest $F$, there is also a \emph{quotient rooted tree} $T/F$. This is a labeled rooted tree with label set $\Gamma$.
Let $U_0$ be the block in $\Gamma$ such that $r$ is in $F(U_0)$. Then $T/F$ is defined on $\Gamma_1 = \Gamma\setminus \{U_0\}$ as follows.
The value of $T/F$ on block $B$  is  obtained by finding the root $j$ of rooted tree $F[B]$ and taking the value to be the
block $B'$ containing $T(j)$.

There is another representation of the quotient rooted tree $T/F$ that may be easier to picture. This is as a labeled rooted tree with label set $R$,
where $R$ is the set of roots in the forest $F$.
The value of $T/F$ on vertex $j$ in $R\setminus \{r\}$  is obtained by finding $T(j)$ and
the block $B'$ that contains it, and taking the value to be the root of rooted tree $F[B']$.

\begin{remark}
The subtree $T_0$ and difference forest $T\setminus T_0$ construction used for the Butcher group
is a special case. The subforest $F$  is $T_0$ together with $F_1 = T \setminus T_0$, and in this case
the quotient rooted tree consists only of a root and immediate successors.
\end{remark}

\begin{example}
 Figure~11
 gives an example of the forests associated with a given rooted tree. The number of rooted trees in the forest (and the number of roots
 of these rooted trees) are 1, 2, 2, 2, 3, 3, 3, 4.
 Figure~12 gives an example
of quotient rooted trees of a given rooted tree on 4 vertices. These correspond to the subforests in Figure~11. The roots of the rooted trees in the forest
give the vertices of the quotient rooted tree. The number of vertices in the quotient rooted trees are 1, 2, 2, 2, 3, 3, 3, 4.
\end{example}

% Figure 11 forests of rooted trees
\begin{figure}[tb]
\begin{picture}(\bggl,\fourgl)(\neggl,0)
\put(0,\twogl){
\begin{picture}(\gl,\gl)
\put(0,0){\line(1,0){\gl}}
\put(0,0){\line(0,1){\gl}}
%\put(0,0){\line(1,1){\gl}}
%\put(\gl,\gl){\line(-1,0){\gl}}
\put(\gl,\gl){\line(0,-1){\gl}}
%\put(0,\gl){\line(1,-1){\gl}}
%\put(-8,0){1}
%\put(45,0){2}
%\put(-8,\gl){3}
%\put(45,\gl){4}
\put(0,0){\circle*{6}}
\put(0,\gl){\circle*{2}}
\put(\gl,0){\circle*{2}}
\put(\gl,\gl){\circle*{2}}
%\put(0,0){\circle*{4}}
%\put(\neghalfgl,\halfgl){4}
\end{picture}}
\put(\twogl,\twogl){
\begin{picture}(\gl,\gl)
%\put(0,0){\line(1,0){\gl}}
\put(0,0){\line(0,1){\gl}}
%\put(0,0){\line(1,1){\gl}}
%\put(\gl,\gl){\line(-1,0){\gl}}
\put(\gl,\gl){\line(0,-1){\gl}}
%\put(0,\gl){\line(1,-1){\gl}}
%\put(-8,0){1}
%\put(45,0){2}
%%\put(-8,\gl){3}
%\put(45,\gl){4}
\put(0,0){\circle*{6}}
\put(0,\gl){\circle*{2}}
\put(\gl,0){\circle*{4}}
\put(\gl,\gl){\circle*{2}}
%\put(0,0){\circle*{4}}
%\put(\neghalfgl,\halfgl){4}
\end{picture}}
\put(\fourgl,\twogl){
\begin{picture}(\gl,\gl)
\put(0,0){\line(1,0){\gl}}
%\put(0,0){\line(0,1){\gl}}
%\put(0,0){\line(1,1){\gl}}
%\put(\gl,\gl){\line(-1,0){\gl}}
\put(\gl,\gl){\line(0,-1){\gl}}
%\put(0,\gl){\line(1,-1){\gl}}
%\put(-8,0){1}
%\put(45,0){2}
%\put(-8,\gl){3}
%\put(45,\gl){4}
\put(0,0){\circle*{6}}
\put(0,\gl){\circle*{4}}
\put(\gl,0){\circle*{2}}
\put(\gl,\gl){\circle*{2}}
%\put(0,0){\circle*{4}}
%\put(\neghalfgl,\halfgl){4}
\end{picture}}
\put(\sixgl,\twogl){
\begin{picture}(\gl,\gl)
\put(0,0){\line(1,0){\gl}}
\put(0,0){\line(0,1){\gl}}
%\put(0,0){\line(1,1){\gl}}
%\put(\gl,\gl){\line(-1,0){\gl}}
%\put(\gl,\gl){\line(0,-1){\gl}}
%\put(0,\gl){\line(1,-1){\gl}}
%\put(-8,0){1}
%\put(45,0){2}
%\put(-8,\gl){3}
%\put(45,\gl){4}
\put(0,0){\circle*{6}}
\put(0,\gl){\circle*{2}}
\put(\gl,0){\circle*{2}}
\put(\gl,\gl){\circle*{4}}
%\put(0,0){\circle*{4}}
%\put(\neghalfgl,\halfgl){4}
\end{picture}}
\put(0,0){
\begin{picture}(\gl,\gl)
%\put(0,0){\line(1,0){\gl}}
%\put(0,0){\line(0,1){\gl}}
%\put(0,0){\line(1,1){\gl}}
%\put(\gl,\gl){\line(-1,0){\gl}}
\put(\gl,\gl){\line(0,-1){\gl}}
%\put(0,\gl){\line(1,-1){\gl}}
%\put(-8,0){1}
%\put(45,0){2}
%\put(-8,\gl){3}
%\put(45,\gl){4}
\put(0,0){\circle*{6}}
\put(0,\gl){\circle*{4}}
\put(\gl,0){\circle*{4}}
\put(\gl,\gl){\circle*{2}}
%\put(0,0){\circle*{4}}
%\put(\neghalfgl,\halfgl){4}
\end{picture}}
\put(\twogl,0){
\begin{picture}(\gl,\gl)
%\put(0,0){\line(1,0){\gl}}
\put(0,0){\line(0,1){\gl}}
%\put(0,0){\line(1,1){\gl}}
%\put(\gl,\gl){\line(-1,0){\gl}}
%\put(\gl,\gl){\line(0,-1){\gl}}
%\put(0,\gl){\line(1,-1){\gl}}
%\put(-8,0){1}
%\put(45,0){2}
%\put(-8,\gl){3}
%\put(45,\gl){4}
\put(0,0){\circle*{6}}
\put(0,\gl){\circle*{2}}
\put(\gl,0){\circle*{4}}
\put(\gl,\gl){\circle*{4}}
%\put(0,0){\circle*{4}}
%\put(\neghalfgl,\halfgl){4}
\end{picture}}
\put(\fourgl,0){
\begin{picture}(\gl,\gl)
\put(0,0){\line(1,0){\gl}}
%\put(0,0){\line(0,1){\gl}}
%\put(0,0){\line(1,1){\gl}}
%\put(\gl,\gl){\line(-1,0){\gl}}
%\put(\gl,\gl){\line(0,-1){\gl}}
%\put(0,\gl){\line(1,-1){\gl}}
%\put(-8,0){1}
%\put(45,0){2}
%\put(-8,\gl){3}
%\put(45,\gl){4}
\put(0,0){\circle*{6}}
\put(0,\gl){\circle*{4}}
\put(\gl,0){\circle*{2}}
\put(\gl,\gl){\circle*{4}}
%\put(0,0){\circle*{4}}
%\put(\neghalfgl,\halfgl){4}
\end{picture}}
\put(\sixgl,0){
\begin{picture}(\gl,\gl)
%\put(0,0){\line(1,0){\gl}}
%\put(0,0){\line(0,1){\gl}}
%\put(0,0){\line(1,1){\gl}}
%\put(\gl,\gl){\line(-1,0){\gl}}
%\put(\gl,\gl){\line(0,-1){\gl}}
%\put(0,\gl){\line(1,-1){\gl}}
%\put(-8,0){1}
%\put(45,0){2}
%\put(-8,\gl){3}
%\put(45,\gl){4}
\put(0,0){\circle*{6}}
\put(0,\gl){\circle*{4}}
\put(\gl,0){\circle*{4}}
\put(\gl,\gl){\circle*{4}}
%\put(0,0){\circle*{4}}
%\put(\neghalfgl,\halfgl){4}
\end{picture}}
\end{picture}
\caption{Subforests}
\end{figure}
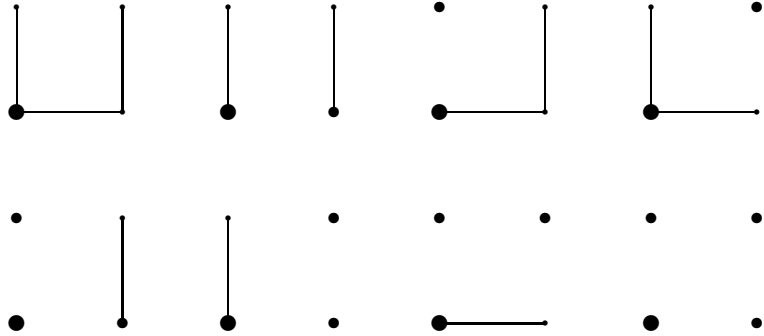

% Figure 12 quotient rooted trees on root sets
\begin{figure}[tb]
\begin{picture}(\bggl,\fourgl)(\neggl,0)
\put(0,\twogl){
\begin{picture}(\gl,\gl)
%\put(0,0){\line(1,0){\gl}}
%\put(0,0){\line(0,1){\gl}}
%\put(0,0){\line(1,1){\gl}}
%\put(\gl,\gl){\line(-1,0){\gl}}
%\put(\gl,\gl){\line(0,-1){\gl}}
%\put(0,\gl){\line(1,-1){\gl}}
%\put(-8,0){1}
%\put(45,0){2}
%\put(-8,\gl){3}
%\put(45,\gl){4}
\put(0,0){\circle*{6}}
%\put(0,\gl){\circle*{2}}
%\put(\gl,0){\circle*{2}}
%\put(\gl,\gl){\circle*{2}}
%\put(0,0){\circle*{4}}
%\put(\neghalfgl,\halfgl){4}
\end{picture}}
\put(\twogl,\twogl){
\begin{picture}(\gl,\gl)
\put(0,0){\line(1,0){\gl}}
%\put(0,0){\line(0,1){\gl}}
%\put(0,0){\line(1,1){\gl}}
%\put(\gl,\gl){\line(-1,0){\gl}}
%\put(\gl,\gl){\line(0,-1){\gl}}
%\put(0,\gl){\line(1,-1){\gl}}
%\put(-8,0){1}
%\put(45,0){2}
%%\put(-8,\gl){3}
%\put(45,\gl){4}
\put(0,0){\circle*{6}}
%\put(0,\gl){\circle*{2}}
\put(\gl,0){\circle*{4}}
%\put(\gl,\gl){\circle*{2}}
%\put(0,0){\circle*{4}}
%\put(\neghalfgl,\halfgl){4}
\end{picture}}
\put(\fourgl,\twogl){
\begin{picture}(\gl,\gl)
%\put(0,0){\line(1,0){\gl}}
\put(0,0){\line(0,1){\gl}}
%\put(0,0){\line(1,1){\gl}}
%\put(\gl,\gl){\line(-1,0){\gl}}
%\put(\gl,\gl){\line(0,-1){\gl}}
%\put(0,\gl){\line(1,-1){\gl}}
%\put(-8,0){1}
%\put(45,0){2}
%\put(-8,\gl){3}
%\put(45,\gl){4}
\put(0,0){\circle*{6}}
\put(0,\gl){\circle*{4}}
%\put(\gl,0){\circle*{2}}
%\put(\gl,\gl){\circle*{2}}
%\put(0,0){\circle*{4}}
%\put(\neghalfgl,\halfgl){4}
\end{picture}}
\put(\sixgl,\twogl){
\begin{picture}(\gl,\gl)
%\put(0,0){\line(1,0){\gl}}
%\put(0,0){\line(0,1){\gl}}
\put(0,0){\line(1,1){\gl}}
%\put(\gl,\gl){\line(-1,0){\gl}}
%\put(\gl,\gl){\line(0,-1){\gl}}
%\put(0,\gl){\line(1,-1){\gl}}
%\put(-8,0){1}
%\put(45,0){2}
%\put(-8,\gl){3}
%\put(45,\gl){4}
\put(0,0){\circle*{6}}
%\put(0,\gl){\circle*{2}}
%\put(\gl,0){\circle*{2}}
\put(\gl,\gl){\circle*{4}}
%\put(0,0){\circle*{4}}
%\put(\neghalfgl,\halfgl){4}
\end{picture}}
\put(0,0){
\begin{picture}(\gl,\gl)
\put(0,0){\line(1,0){\gl}}
\put(0,0){\line(0,1){\gl}}
%\put(0,0){\line(1,1){\gl}}
%\put(\gl,\gl){\line(-1,0){\gl}}
%\put(\gl,\gl){\line(0,-1){\gl}}
%\put(0,\gl){\line(1,-1){\gl}}
%\put(-8,0){1}
%\put(45,0){2}
%\put(-8,\gl){3}
%\put(45,\gl){4}
\put(0,0){\circle*{6}}
\put(0,\gl){\circle*{4}}
\put(\gl,0){\circle*{4}}
%\put(\gl,\gl){\circle*{2}}
%\put(0,0){\circle*{4}}
%\put(\neghalfgl,\halfgl){4}
\end{picture}}
\put(\twogl,0){
\begin{picture}(\gl,\gl)
\put(0,0){\line(1,0){\gl}}
%\put(0,0){\line(0,1){\gl}}
%\put(0,0){\line(1,1){\gl}}
%\put(\gl,\gl){\line(-1,0){\gl}}
\put(\gl,\gl){\line(0,-1){\gl}}
%\put(0,\gl){\line(1,-1){\gl}}
%\put(-8,0){1}
%\put(45,0){2}
%\put(-8,\gl){3}
%\put(45,\gl){4}
\put(0,0){\circle*{6}}
%\put(0,\gl){\circle*{2}}
\put(\gl,0){\circle*{4}}
\put(\gl,\gl){\circle*{4}}
%\put(0,0){\circle*{4}}
%\put(\neghalfgl,\halfgl){4}
\end{picture}}
\put(\fourgl,0){
\begin{picture}(\gl,\gl)
%\put(0,0){\line(1,0){\gl}}
\put(0,0){\line(0,1){\gl}}
\put(0,0){\line(1,1){\gl}}
%\put(\gl,\gl){\line(-1,0){\gl}}
%\put(\gl,\gl){\line(0,-1){\gl}}
%\put(0,\gl){\line(1,-1){\gl}}
%\put(-8,0){1}
%\put(45,0){2}
%\put(-8,\gl){3}
%\put(45,\gl){4}
\put(0,0){\circle*{6}}
\put(0,\gl){\circle*{4}}
%\put(\gl,0){\circle*{2}}
\put(\gl,\gl){\circle*{4}}
%\put(0,0){\circle*{4}}
%\put(\neghalfgl,\halfgl){4}
\end{picture}}
\put(\sixgl,0){
\begin{picture}(\gl,\gl)
\put(0,0){\line(1,0){\gl}}
\put(0,0){\line(0,1){\gl}}
%\put(0,0){\line(1,1){\gl}}
%\put(\gl,\gl){\line(-1,0){\gl}}
\put(\gl,\gl){\line(0,-1){\gl}}
%\put(0,\gl){\line(1,-1){\gl}}
%\put(-8,0){1}
%\put(45,0){2}
%\put(-8,\gl){3}
%\put(45,\gl){4}
\put(0,0){\circle*{6}}
\put(0,\gl){\circle*{4}}
\put(\gl,0){\circle*{4}}
\put(\gl,\gl){\circle*{4}}
%\put(0,0){\circle*{4}}
%\put(\neghalfgl,\halfgl){4}
\end{picture}}
\end{picture}
\caption{Quotient rooted trees}
\end{figure}
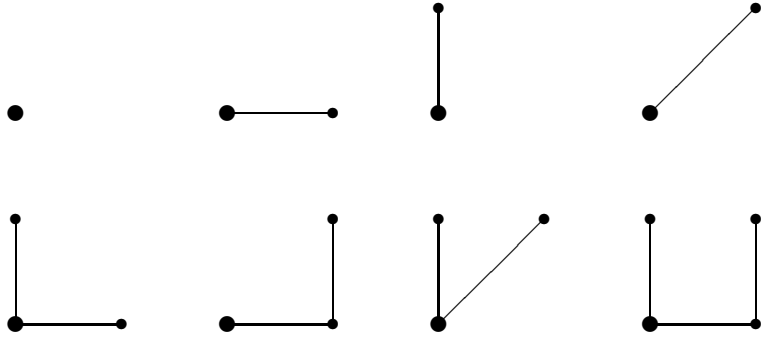

The pair consisting of
the subforest $F$ and the quotient rooted tree $T/F$ do not completely determine the original rooted tree $T$. The quotient rooted tree
assigns to each block in $\Gamma_1$ another block in $\Gamma$. It is also necessary to specify a point in that block.

\begin{proposition}
For each label set $U$ there is a one-to-one correspondence between rooted tree, subforest pairs $T, F$ and triples $F$, $\hat T$, $\phi$, where
$F$ is a forest of rooted trees with set partition $\Gamma$, $\hat T$ is a rooted tree on vertex set $\Gamma$ with root $U_0$,
and $\phi$ is a function defined on $\Gamma_1 = \Gamma \setminus \{U_0\}$ such that for every $B$ in $\Gamma_1$ the value $\phi(B)$ is in $\hat T(B)$.
\end{proposition}

 In the proposition the function $\phi$ sends a block to a vertex in another target block.
There is a parametrization of such functions by  target blocks.
If $B'$ is a block in $\Gamma$, then define $\Phi(B')$ to be the restriction of $\phi$ to
$T^{-1}(B')$. Thus if $\hat T(B) = B'$, then $\Phi(B')$ applied to $B$ is a vertex in $B'$. Conversely, suppose that there is a function $\Phi$
that maps each $B'$ in $\Gamma$ to a function $\Phi(B')$ from $\hat T^{-1}(B')$ to $B'$. Then there is a corresponding $\phi$
given on $B$ with $\hat T(B) = B'$ by $\phi( B) = \Phi(B')(B)$.

Let $\tilde A_\emptyset$ be the set of unlabeled rooted trees augmented with the empty object. There is a multiplication defined on certain
elements of $\Rb^{\tilde \Ac_\emptyset}$.
This is the \emph{quotient rooted tree convolution} $c = a \star b$, defined whenever $a(\emptyset) = 0$, such that
\begin{equation}
c(T) =   \sum_{F \sqsubseteq T}  b( T/F) a^\times( F)  .
\end{equation}
The sum is over subforests $F$ of $T$. The contribution of a forest is $a^\times(F) = \prod_{T' \in F} a(T')$. For the empty
forest this product is 1. For the special case of the empty set object $c(\emptyset) = b(\emptyset)a^\times(\emptyset) = b(\emptyset)$.

If $R = \{r\}$, then  the forest $F$ has only one rooted tree $T$, and $T/F$ is a rooted tree on a one point vertex set. So the contribution to the sum
is $a(T) b(\bullet) $.
 When $R$ is the entire vertex set, then  $F$ is the discrete forest, and $T/F = T$.
 The contribution to the sum is $a(\bullet)^{|T|}b(T) $. As a special case, $c(\bullet) = a(\bullet) b(\bullet)$.

 If the multiplication is restricted to  $a(\emptyset) = b(\emptyset) = 0$, then $G_S = \Rb^{\tilde \Ac}$
  becomes an algebraic system closed under multiplication. The $S$ stands for substitution.
 The identity for quotient rooted tree convolution is $\delta_\bullet$, which has
 coefficient 1 for a one-vertex rooted tree and 0 for all other rooted trees.
 If the multiplication is also restricted to  $a(\bullet)   \neq 0,  b(\bullet) \neq  0$, the resulting system $G_S^\star $
is a group.  If the multiplication is further restricted to  $a(\bullet) = b(\bullet) = 1$, then this defines a subgroup $G_S^1$.
The group $G_S^1$ is the character group of the  rooted tree Hopf algebra of Calaque,  Ebrahimi-Fard, and Manchon.

%It is possible to extend the theory to the case when the only requirement is that $a(\emptyset) = 0$. If $T$ is
%not the empty rooted tree, then there is no change.
%When $T$ is the empty rooted tree, then $R = \emptyset$, $F$ is the empty forest (contributing 1), and $T/F$ is
% the empty rooted tree. So $c(\emptyset) = b(\emptyset)$.

There is also a functional representation of quotient rooted tree convolution. This deals with formal functions of the
form $g \mapsto \alpha(g) = \sum_{n=0}^\infty (1/n!) a_n g^n$. Such a function is denoted $\alpha$ to indicate that it
depends only on the coefficients $a_n$ and not on the input $g$.
The power series depends on the choice of $\alpha$ and is of the form
\begin{equation}
f^c(t, \alpha, g) = \sum_{n=0}^\infty \sum_{T \in \Ac[U_n] } \frac{t^n}{n!} c(T) \prod_{\ell \in U_n} \alpha^{(|T^{-1}(\ell)|)}(g).
\end{equation}
The zero order term is $c(0) g$.
This leads to the following theorem \cite{ChartierHairerVilmart}.

\begin{theorem}[Substitution] Suppose that $a(\emptyset) = 0$ and $c = a \star b$ is the quotient rooted tree convolution.  Then
\begin{equation}
f^c( t , \beta , g) = f^b( 1,  f^a(t ,\beta, \cdot), g) .
\end{equation}
\end{theorem}

\begin{proof} The proof here uses
 the calculus formulas in the appendix. Let $D_t$ be the partial derivative with respect to $t$. It is sufficient to show
that applying $D_t^n$ and then setting $t$ to zero gives the same result for both sides of the equation.
For the left hand side this is
\begin{equation}
D_t^{|U|} f^c( t,\beta, g) |_{t = 0} =
\sum_{T \in \Ac[U] }  c(T) \prod_{\ell \in U} D_g^{|T^{-1}(\ell)|} \beta(g),
\end{equation}
where $U$ is a label set with $n$ points.

\noindent\textbf{Product rule} The computation for the right hand side  begins with
\begin{equation}
f^b( 1,  f^a (t, \beta, \cdot), g) = \sum_{\hat T}\frac{1}{|\hat T|!}   b(\hat T) \prod_{i \in [\hat T]} D_g^{|{\hat T}^{-1}(i)|}  f^a(t,\beta, g)
\end{equation}
By the product rule
\begin{equation}
D_t^{|U|}  \prod_{i \in [\hat T]} D_g^{|{\hat T}^{-1}(i)|} f^a (t ,\beta, g)   =
\sum_{\psi: U \to [\hat T] }
\prod_{i \in [\hat T]}   D_t^{|\psi^{-1}(i)|} D_g^{|{\hat T}^{-1}(i)|} f^a(t,\beta,g)  .
\end{equation}
Set $t = 0$. Since $f^a(0,\beta,g) = a(\emptyset) g = 0$, the contributions from $i$ with $\psi^{-1}(i) = \emptyset $ are zero. So $\psi: U \to [\hat T]$ induces a set partition $\Gamma$ of $U$
and a bijection from $\Gamma$ to $[\hat T]$.
There are $|\hat T|!$ such bijections.
This leads to
\begin{equation}
D_t^{|U|} f^b(1,  f^a (t,\beta, \cdot), g) |_{t=0} = \sum_{\Gamma \in \setpartition[U]} \sum_{\hat T \in \Ac[\Gamma]}  c(\Gamma, \hat T) ,
\end{equation}
where
\begin{equation}
c(\Gamma, \hat T) =
 \prod_{B \in \Gamma} b(\hat T)
D_g^{{|\hat T}^{-1}(B)|} D_t^{|B|}  f^a(t, \beta, g)|_{t=0} .
\end{equation}

\noindent\textbf{Distributive law}
The next stage is to insert
\begin{equation}
D_t^{|B|} f^a( t,\beta, g) |_{t = 0} =
\sum_{H \in \Ac[B] }  a(H) \prod_{j \in B} D_g^{|H^{-1}(j)|} \beta(g)
\end{equation}
and use the distributive law. This gives a forest sum
\begin{equation}
D_t^{|U|} f^b(1,  f^a (t,\beta, \cdot), g) |_{t=0} = \sum_{\Gamma \in \setpartition[U]} \sum_{\hat T \in \Ac[\Gamma]} \sum_F
c(\Gamma, \hat T, F) ,
\end{equation}
where
\begin{equation}
c(\Gamma, \hat T, F) = \prod_{B \in \Gamma}
 b(\hat T)  a(F(B))
D_g^{|{\hat T}^{-1}(B)|} \prod_{j \in B} D_g^{|F(B)^{-1}(j)|} \beta(g) .
\end{equation}

\noindent\textbf{Product rule}
The product rule for differentiation produces
\begin{equation}
c(\Gamma, \hat T, F)
= \prod_{B \in \Gamma}
b(\hat T) a(F(B)) \sum_{\phi: {\hat T}^{-1}(B)) \to B} \prod_{j \in B} D_g^{|\phi^{-1}(j)|} D_g^{|F(B)^{-1}(j|} \beta(g) .
\end{equation}

\noindent\textbf{Distributive law}
The distributive law gives
\begin{equation}
D_t^{|U|} f^b(1,  f^a (t,\beta, \cdot), g) |_{t=0} = \sum_{\Gamma \in \setpartition[U]} \sum_{\hat T \in \Ac[\Gamma]} \sum_F
\sum_\Phi
c(\Gamma, \hat T, F, \Phi) ,
\end{equation}
where
\begin{equation}
c(\Gamma, \hat T, F, \Phi) = \prod_{B \in \Gamma} \prod_{j \in B}
 b(\hat T)  a(F(B))
 D_g^{|\Phi(B)^{-1}(j)|} D_g^{|F(B)^{-1}(j|)} \beta(g) .
\end{equation}

The sum is over set partitions $\Gamma$ of $U$ and over forest functions $F$ that send block $B$ to rooted tree $F(B)$ on vertex set $B$.
It is also over rooted trees $\hat T$ on vertex set $\Gamma$. Finally, it is over functions $\Phi$ that
send each block $B$ in $\Gamma$ to a function $\Phi(B)$ that takes each $\hat T$ preimage block $B'$ and sends it to a vertex in $B$.

\noindent\textbf{Construction of rooted tree and subforest}
These data  determine a rooted tree on $U$ that is made from the rooted trees $F(B)$ internal to the blocks $B$ and from the rooted tree $\hat T$ and the function $\Phi$. If $\hat T(B') = B$, then there is an edge from the root of the rooted tree on $B'$ to the $\Phi(B)(B')$ in $B$.

The corresponding contribution involves the coefficients $b(\hat T)$
and $a(F(B))$ and a contribution from the rooted tree. At a given vertex $j$ this involves a derivative of $\beta(g)$ of an order equal to the
total number of edges impinging on this vertex, both from within the block and from the other blocks.

Giving these data is the same as giving the pair $T \in \Ac[U]$ together with subforest $F$. So the final expression
is
\begin{equation}
D_t^{|U|} f^b(  f^a (t ,\beta, \cdot), g)|_{t=0} = \sum_{T \in \Ac[U]} \sum_{F \sqsubseteq T}  a^\times(F) b(T/F)  f_T(g).
\end{equation}
This establishes the result.
\end{proof}

The authors \cite{ChartierHairerVilmart} give two applications of this result. For both the idea is to consider the coefficients
$e(\tau) = 1/\tau!$ that give the exact solution of an ordinary differential equation $dx/dt = \beta(x)$. In backward error
analysis the idea is to take $c$ corresponding to some numerical method and solve $c = a \star e$ for $a$. This produces a
modified differential equation that agrees with the numerical method. In the application to modified integrators, start with
a numerical method given by $b$ and solve $e = a \star b$ for $a$. This produces a modified numerical method that agrees
with the solution of the differential equation.

\iffalse
The authors \cite{ChartierHairerVilmart} give two applications of this result. The first is to backward error analysis.
Let $b(\tau) = 1/\tau!$, so $x = f^b(t,\beta, g)$ be the exact
solution of the ordinary differential equation $dx/dt = \beta(x)$, with $x = g$ at $t = 0$. Let $f^c(h,\beta, g)$ be a numerical method for
solving the same equation with time step $h$.  Since $G_S^1$ is a group, there is an $f^a(h,\beta, \cdot)$
with $f^c( h , \beta , g) = f^b( 1,  f^a(h ,\beta, \cdot), g) $. This says that the result of the  numerical method agrees with
the exact solution of the modified differential equation  $dx/dt = f^a(h, \beta, x) $ with $x = g$ at $t=0$, evaluated at $t = 1$.

The second application is to modified integrators. Let $c(\tau) = 1/\tau!$, so $x = f^c(t,\beta, g)$ is the exact
solution of the ordinary differential equation $dx/dt = \beta(x)$, with $x = g$ at $t = 0$. Let $f^b(h,\beta, g)$ be a numerical method for
solving the same equation with time step $h$.  Since $G_S^1$ is a group, there is a  $f^a(t,\beta, \cdot)$
with $f^c( t , \beta , g) = f^b( 1,  f^a(t ,\beta, \cdot), g) $. This says that the  solution of the  differential equation
agrees with the result of the modified numerical method in which $\beta$ is replaced by $f^a(t,\beta, \cdot)$, evaluated at $h=1$.
\fi

\section*{Appendix: Algebra and calculus in combinatorics}

Here are some basic results from algebra and calculus in forms that are useful
for combinatorics. These are stated in the setting of functions of one variable. There
are even more illuminating multi-variable results, but they are not needed in the present
exposition.

\noindent\textbf{The distributive law}
A version of the  distributive law of algebra is the following. Suppose that $B$ is a set and for each $b \in B$ there is
a corresponding set $F_b$.  Then the product over $b \in B$ of sums indexed by $F_b$ is a sum of products:
\begin{equation}
 \prod_{b \in B} \sum_{ t \in F_b} a_b(t)    = \sum_{s }  \prod_{b \in B} a_b(s(b)).
\end{equation}
The sum on the right is over  all functions $s: B \to \bigcup_b F_b$ such that
for each $b$ the value $s(b) \in F_b$. The set of all such functions is the product space $\prod_b F_b$. In the special case when the $F_b=F$ are all the same, the sum is over all functions $s: B \to F$. In this case the set of all such functions is the Cartesian power space $F^B$.

\noindent\textbf{The product rule}
A version of the  product rule for differentiation is the following. Let $U$ be a set with $|U|$ elements. Then the $|U|$ order
derivative of a product function is given by
\begin{equation}
D^{|U|} \prod_{b \in B} F_b  = \sum_{\phi \in B^U} \prod_{b \in B} D^{ |\phi^{-1}(b)|}F_b.
\end{equation}
Here $B^U$ consists of all functions $\phi: U \to B$. Sometimes a function $\phi$ is described
by its inverse images $U_b = \phi^{-1}(b)$, so one can think of this as a sum over the corresponding maps $b \mapsto U_b$.

\noindent\textbf{The chain rule}
A version of the chain rule is the following. Let $U$ be a set with $|U|$ elements. Then the $|U|$ order
derivative of a composite function is given by
\begin{equation}
D^{|U|}(F \circ G)  = \sum_{\Gamma \in \setpartition[U]} (D^{|\Gamma|}F)\circ G \cdot \prod_{B \in \Gamma} D^{|B|} G.
\end{equation}
Here $\setpartition[U]$ consists of all set partitions of $U$ into disjoint non-empty subsets with union $U$.

\end{document}